\documentclass[envcountsame,envcountresetchap]{svmult}

\usepackage{makeidx}
\usepackage{graphicx}
\usepackage{multicol}
\usepackage[bottom]{footmisc}

\usepackage{amssymb}
\usepackage{amsmath}
\usepackage{dsfont}
\RequirePackage[T1]{fontenc}

\smartqed

\input xy
\xyoption{all}

\begin{document}
%
%
\newtheorem{Definition}{Definition}[section]
\newtheorem{Conclusion}{Conclusion}[section]
\newtheorem{Construction}{Construction}[section]
\newtheorem{Example}[Definition]{Example}
\newtheorem{Examples}[Definition]{Examples}
\newtheorem{Remark}[Definition]{Remark}
\newtheorem{Remarks}[Definition]{Remarks}
\newtheorem{Caution}[Definition]{Caution}
\newtheorem{Conjecture}[Definition]{Conjecture}
\newtheorem{Question}[Definition]{Question}
\newtheorem{Questions}[Definition]{Questions}
\newtheorem{Theorem}[Definition]{Theorem}
\newtheorem{Proposition}[Definition]{Proposition}
\newtheorem{Lemma}[Definition]{Lemma}
\newtheorem{Corollary}[Definition]{Corollary}
\newtheorem{Fact}[Definition]{Fact}
\newtheorem{Facts}[Definition]{Facts}

\newcommand{\prf}{\par\noindent{\sc Proof.}\quad}
\newcommand{\blowup}{\rule[-3mm]{0mm}{0mm}}
\newcommand{\Aff}{{\mathds{A}}}
\newcommand{\BB}{{\mathds{B}}}
\newcommand{\CC}{{\mathds{C}}}
\newcommand{\FF}{{\mathds{F}}}
\newcommand{\GG}{{\mathds{G}}}
\newcommand{\HH}{{\mathds{H}}}
\newcommand{\NN}{{\mathds{N}}}
\newcommand{\ZZ}{{\mathds{Z}}}
\newcommand{\PP}{{\mathds{P}}}
\newcommand{\QQ}{{\mathds{Q}}}
\newcommand{\RR}{{\mathds{R}}}
\newcommand{\Liea}{{\mathfrak a}}
\newcommand{\Lieb}{{\mathfrak b}}
\newcommand{\Lieg}{{\mathfrak g}}
\newcommand{\Liem}{{\mathfrak m}}
\newcommand{\ideala}{{\mathfrak a}}
\newcommand{\idealb}{{\mathfrak b}}
\newcommand{\idealg}{{\mathfrak g}}
\newcommand{\idealm}{{\mathfrak m}}
\newcommand{\idealp}{{\mathfrak p}}
\newcommand{\idealq}{{\mathfrak q}}
\newcommand{\idealI}{{\cal I}}
\newcommand{\lin}{\sim}
\newcommand{\num}{\equiv}
\newcommand{\dual}{\ast}
\newcommand{\iso}{\cong}
\newcommand{\homeo}{\approx}
\newcommand{\mm}{{\mathfrak m}}
\newcommand{\pp}{{\mathfrak p}}
\newcommand{\qq}{{\mathfrak q}}
\newcommand{\rr}{{\mathfrak r}}
\newcommand{\pP}{{\mathfrak P}}
\newcommand{\qQ}{{\mathfrak Q}}
\newcommand{\rR}{{\mathfrak R}}
%
%
\newcommand{\dq}{{``}}
\newcommand{\OO}{{\cal O}}
\newcommand{\into}{{\hookrightarrow}}
\newcommand{\onto}{{\twoheadrightarrow}}
\newcommand{\Spec}{{\rm Spec}\:}
\newcommand{\BigSpec}{{\rm\bf Spec}\:}
\newcommand{\Proj}{{\rm Proj}\:}
\newcommand{\Pic}{{\rm Pic }}
\newcommand{\Br}{{\rm Br}}
\newcommand{\NS}{{\rm NS}}
\newcommand{\chit}{\chi_{\rm top}}
\newcommand{\KanDiv}{{\cal K}}
\newcommand{\perdef}{{\stackrel{{\rm def}}{=}}}
\newcommand{\Cycl}[1]{{\ZZ/{#1}\ZZ}}
\newcommand{\Sym}{{\mathfrak S}}
\newcommand{\Xcan}{X_{{\rm can}}}
\newcommand{\Ycan}{Y_{{\rm can}}}
\newcommand{\ab}{{\rm ab}}
\newcommand{\Aut}{{\rm Aut}}
\newcommand{\Hom}{{\rm Hom}}
\newcommand{\Supp}{{\rm Supp}}
\newcommand{\ord}{{\rm ord}}
\newcommand{\divisor}{{\rm div}}
\newcommand{\Alb}{{\rm Alb}}
\newcommand{\Jac}{{\rm Jac}}
\newcommand{\Ig}{{\rm Ig}}
\newcommand{\piet}{{\pi_1^{\rm \acute{e}t}}}
\newcommand{\Het}[1]{{H_{\rm \acute{e}t}^{{#1}}}}
\newcommand{\Hcris}[1]{{H_{\rm cris}^{{#1}}}}
\newcommand{\HdR}[1]{{H_{\rm dR}^{{#1}}}}
\newcommand{\hdR}[1]{{h_{\rm dR}^{{#1}}}}

\title*{Algebraic Surfaces in Positive Characteristic}
\titlerunning{Algebraic Surfaces}

\author{Christian Liedtke}
\authorrunning{C.~Liedtke}
\institute{Mathematisches Institut\\ 
Endenicher Allee 60\\ 
D-53115 Bonn, Germany\\
\texttt{liedtke@math.uni-bonn.de}}

\setcounter{tocdepth}{1}
\tableofcontents

\maketitle

\begin{abstract}
 These notes are an introduction to and an overview over the theory
 of algebraic surfaces over algebraically closed fields of positive characteristic.
 After a little bit of background in characteristic-$p$-geometry,
 we sketch the Kodaira--Enriques classification.
 Next, we turn to more special characteristic-$p$ topics, 
 and end with lifting results, as well as applications to geometry in characteristic zero.
 We assume that the reader has a background in complex geometry and has
 seen the Kodaira--Enriques classification of complex surfaces
 before.
\end{abstract}


\


\section{Introduction}

These notes grew out of a series of lectures given at Sogang University,
Seoul in October 2009.
They were meant for complex geometers, who are  not familiar with
characteristic-$p$-geometry, but who would like to see similarities,
as well as differences to complex geometry.
More precisely, these notes are on algebraic surfaces
in positive characteristic, and assume familiarity with the complex
side of this theory, say,
on the level of Beauville's book \cite{beauville book}.

Roughly speaking, the theory of curves in characteristic
zero and $p>0$ look very similar, and many fundamental
results were already classically known to hold in arbitrary characteristic.
Also, curves lift from characteristic $p$ over
the Witt ring to characteristic zero, which implies that
many ``characteristic-$p$-pathologies'' cannot happen.
Abelian varieties admit at least formal lifts over the Witt ring,
and we refer to Section \ref{sec:Witt} for details and implications.

However, from dimension two on, geometry 
in positive characteristic displays more and more differences to classical 
complex geometry.
In fact, this geometry has long been considered as ``pathological'' and
``exotic'', as even reflected in the titles of a series of articles by Mumford 
(the first one being \cite{Mumford Path1}).
There, the emphasis was more on finding and exploring differences
to the classical theory.
For a short overview over the main new phenomena
for surfaces in positive characteristic, we refer to \cite[Section 15]{isk shaf}.

However, in their three fundamental articles \cite{Mumford}, \cite{bm2} and \cite{bm3},
Bombieri and Mumford established the Kodaira--Enriques classification
of algebraic surfaces in positive characteristic.
Together with Artin's results \cite{Artin Sing1} and \cite{Artin Sing2}
on surface singularities,  especially rational and Du~Val singularities, 
as well as work of Ekedahl \cite{ek} on pluricanonical systems of surfaces
of general type (extending Bombieri's results to characteristic $p$),
this sets the scene in positive characteristic.
It turns out that surface theory in positive characteristic is in many respects 
not so different from characteristic zero, at least, if one takes the right angle.

As over the complex numbers, there is a vast number of examples,
counter-examples, and
(partial) classification results for special classes of surfaces in
positive characteristic.
Unfortunately, it was impossible for me to mention all of them in these introductory notes -- 
for example, I could have written  much more on K3 surfaces, elliptic
surfaces, and (birational) automorphisms of surfaces.
\medskip

These notes are organized as follows:
\begin{trivlist}
 \item[{\bf Preparatory Material}]
 \item[{\em Section \ref{sec:Curves}}] We introduce the
  various Frobenius morphisms,
  and proceed to basic results on algebraic curves.
  Finally, we discuss finite, constant, and
  infinitesimal group schemes, as well as the three group schemes of length $p$.
  \item[{\em Section \ref{sec:Cohomology}}] We recall Hodge-, \'etale 
   and deRham (hyper-) cohomology.
   Next, we discuss Albanese and Picard schemes, non-closed 
   differential forms,
   and their relation to  (non-)degeneracy of the Fr\"olicher 
   spectral sequence from Hodge- to deRham-cohomology.
   Finally, we sketch how crystalline
   cohomology links all the above mentioned cohomology theories.
\medskip

  \item[{\bf Classification of Algebraic Surfaces}] 
  \item[{\em Section \ref{sec:Birational}}] 
   We discuss blow-ups and Castelnuovo's contraction theorem, 
   introduce minimal models,
   and describe the structure of rational and birational maps of surfaces.
   We classify birationally ruled surfaces, and state the rationality theorem
   of Castelnuovo--Zariski.
  \item[{\em Section \ref{sec:Elliptic}}]
   We recall elliptic fibrations, and discuss the phenomenona
   of quasi-elliptic fibrations and wild fibers.
   Then, we state the canonical bundle formula and give the 
   possible degeneration types of fibers in (quasi-)elliptic fibrations.
  \item[{\em Section \ref{sec:Classification}}]
   We sketch the Kodaira--Enriques classification of algebraic surfaces
   according to their Kodaira dimension.
  \item[{\em Section \ref{sec:KodairaZero}}] We discuss the four classes of
   minimal surfaces of Kodaira dimension zero in greater detail.
   We put an emphasis on the non-classical classes of Enriques surfaces
   in characteristic $2$,
   as well as the new classes of quasi-hyperelliptic surfaces in
   characteristics $2$ and $3$.
  \item[{\em Section \ref{sec:GeneralType}}]
   We start with Ekedahl's work on pluricanonical maps of surfaces of general type.
   Then, we continue with what is known about various inequalities (Noether,
   Castelnuovo, Bogomolov--Miyaoka--Yau) in positive characteristic,
   and end with a couple of results on surfaces of general type with small
   invariants.
 \medskip

  \item[{\bf Special Topics in Positive Characteristic}]
  \item[{\em Section \ref{uni and super}}]   
  We study uniruled surfaces that are not birationally ruled,
   and introduce two notions of supersingularity, due to Artin and Shioda.
   Then, we discuss these notions for K3 surfaces.
   Next, we turn to surfaces
   over finite fields, zeta functions, and the Tate conjecture.
  \item[{\em Section \ref{sec:Foliations}}]
   We explain Jacobson's correspondence for purely inseparable
   field extensions. 
   On the geometric level, this corresponds to $p$-closed foliations.
    We give applications to global vector fields on K3 surfaces and 
    end by discussing quotients by infinitesimal group schemes.
 \medskip
 
   \item[{\bf  From Positive Characteristic to Characteristic Zero}]
   \item[{\em Section \ref{sec:Witt}}]
    We recall the ring of Witt vectors, discuss various notions
    of what it means to ``lift to characteristic zero'', and discuss,
    what each type of liftings implies.
    We end by giving examples and counter-examples.
   \item[{\em Section \ref{sec:rationalcurves}}]
    As an application of characteristic-$p$-geometry, we establish
    infinitely many rational curves on complex projective K3 surfaces
    of odd Picard rank using reduction modulo $p$ and
    special characteristic-$p$ features.
\end{trivlist}

Finally, we advise the reader who is interested in learning
surface theory over algebraically closed 
ground fields of arbitrary characteristic  from scratch (including proofs)
to have a look at Badescu's excellent text book \cite{Badescu}.
From there, the reader can proceed to more advanced topics, including
the original articles by Bombieri and Mumford mentioned above,
as well as the literature given in these notes.
\medskip

\noindent {\bf Acknowledgments.} 
These lecture notes grew out of a lecture series given at Sogang University,
Seoul in October 19-22, 2009.
I thank Yongnam Lee for the invitation to Sogang University and
hospitality.
It was a pleasure visiting him and giving these lectures.
Also, I thank Fabrizio Catanese, Hisanori Ohashi,
Holger Partsch, S\"onke Rollenske, Nguyen Le Dang Thi,
Yuri Tschinkel, Tong Zhang, and the referee for suggestions, remarks, and 
pointing out mistakes in earlier versions.
I thank the referee especially for clarifications and providing me with more references.
I wrote up a first version of these notes at Stanford University and
I thank the department for hospitality.
I gratefully acknowledge funding from DFG under research grants 
LI 1906/1-1 and LI 1906/1-2.

\addtocontents{toc}{{\medskip \bf Preparatory Material}}

\section{Frobenius, curves and group schemes}
\label{sec:Groups}
\label{sec:Curves}

Before dealing with surfaces, we first shortly review
a little bit of background material.
Of course, the omnipresent Frobenius morphism has to
be mentioned first -- in many cases, when a 
characteristic zero argument breaks down in
positive characteristic, inseparable morphisms and
inseparable field extensions are responsible.
The prototype of an inseparable morphism is
the {\em Frobenius morphism}, and in
many situations it also provides the solution 
to a problem.
Next, we give a very short overview over curves 
and group scheme actions.
We have chosen our material in view of what we
need for the classification and description of
surfaces later on.

\subsection{Frobenius}
\label{subsec:frobenius}
Let us recall that a field $k$ of positive characteristic $p$
is called {\em perfect}, if its 
Frobenius morphism $x\mapsto x^p$
is surjective, i.e., if every element in $k$ has a $p$.th root in $k$.
For example, finite fields and algebraically closed fields are perfect.
On the other hand, function fields of varieties in positive
characteristic are almost never perfect.

Let $X$ be an $n$-dimensional variety over a field $k$
with structure morphism $f:X\to\Spec k$.
Then, the morphism that is the identity on points of $X$
and is $x\mapsto x^p$ on the structure sheaf $\OO_X$ is called the
{\em absolute Frobenius morphism} $F_{\rm abs}$ of $X$.

However, the absolute Frobenius morphism is not ``geometric'':
namely, it acts as
$x\mapsto x^p$ on the ground field $k$,
which is non-trivial except for $k=\FF_p$.
To obtain a morphism over $k$, we first form the pull-back
$$
 X^{(p)} \,:=\, X\times_{\Spec k}\Spec k 
 \,\stackrel{ {\rm pr}_2 }{\longrightarrow} \Spec k 
$$
with respect to the structure map $f:X\to\Spec k$ and with respect
to the absolute Frobenius $F_{\rm abs}:\Spec k\to\Spec k$.
This {\em Frobenius pullback}
$f^{(p)}:={\rm pr}_2:X^{(p)}\to \Spec k$ 
is a new variety over $k$. 
If $k$ is perfect then $X$ and $X^{(p)}$ are
abstractly isomorphic as schemes, but not as varieties over $k$,
see below.

Now, by the universal property of pull-backs, we obtain a morphism
of varieties over $k$, the {\em $k$-linear Frobenius morphism}
$F:X\to X^{(p)}$
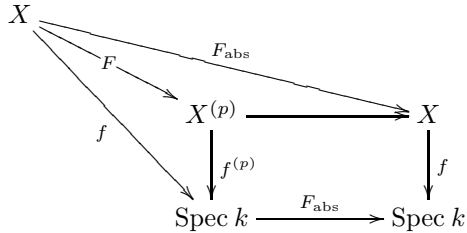
$$
\xymatrixcolsep{4pc}
\xymatrix{
X \ar[ddr]_f \ar[drr]^{F_{\rm abs}}
\ar[dr]|-{F} \\
& X^{(p)} \ar[d]^{f^{(p)}} \ar[r]
& X\ar[d]^{f}\\
&\Spec k \ar[r]^{F_{\rm abs}} & \Spec k}
$$
In more down to earth terms and for affine space this simply means
$$\begin{array}{lccccl}
              &&  k[x_1,...,x_n] & \to & k[x_1,...,x_n] \\
  \mbox{ absolute Frobenius } &:& f(x_1,...,x_n) &\mapsto & \left(f(x_1,...,x_n)\right)^p\\
  \mbox{ $k$-linear Frobenius }&:& c &\mapsto & c & \mbox{ if } c\in k\\
    && x_i &\mapsto & x_i^p
  \end{array}
$$
When dealing with varieties over finite fields there are 
even more Frobenius morphisms: over the field $\FF_q$ with
$q=p^n$ elements one has a Frobenius morphism $F_q:x\mapsto x^q$, 
and for technical reasons sometimes its inverse has to be considered,
see, for example \cite[Appendix C.4]{Hartshorne}.
Depending on author and context all these morphisms and various
base-changes are called ``Frobenius'' and so, a little 
care is needed.

Next, if $X$ is $n$-dimensional over $k$, then the $k$-linear Frobenius
$F:X\to X^{(p)}$ is a finite morphism of degree $p^n$.
Moreover, if $k$ is perfect then, on the level of function fields,
this morphism corresponds to the inclusion
$$
     k(X^{(p)})\,=\,k(X)^p\,\subseteq\, k(X) \,.
$$
Note that $k(X)^p$, the set of $p$.th powers of $k(X)$, is in fact a
field: it is not only closed under multiplication, but also under
addition since $x^p+y^p=(x+y)^p$ in characteristic $p$.
Let us also fix an algebraic closure $\overline{K}$ of $K=k(X)$.
For every integer $i\geq0$ we define
$$
   K^{p^{-i}} \,:=\, \{\, x\in\overline{K}\,|\, x^{p^i}\in K\,\}
$$
and note that these sets are in fact fields.
The field $K^{p^{-i}}$ is a finite and purely inseparable
extension of $K$ of degree $p^{ni}$.
Their union $K^{p^{-\infty}}$ as $i$ tends to infinity is called the 
{\em perfect closure} of $K$ in $\overline{K}$, as it is
the smallest subfield of $\overline{K}$
that is perfect and contains $K$.

\begin{Definition}
 Let $L$ be a finite and purely inseparable extension of $K$.
 The {\em height} of $L$ over $K$ is defined to be the minimal
 $i$ such that $K\subseteq L\subseteq K^{p^{-i}}$.
\end{Definition}

Similarly, if $\varphi:Y\to X$ is a purely inseparable and
generically finite morphism of varieties over a perfect field $k$,
then the {\em height of $\varphi$} is defined to be the
height of the extension of function fields $k(Y)/k(X)$.
For example, the $k$-linear Frobenius morphism is of height one.

For more on inseparable morphisms, we refer to 
Section \ref{sec:Foliations}.

\subsection{Curves}
 \label{subsec:curves}
 Most of the results of this section are classical, and we refer
to \cite[Chapter IV]{Hartshorne} or \cite{Liu} for details,
specialized topics and further references.
Let $C$ be a smooth projective curve over an algebraically
closed field of characteristic $p\geq0$.
Then its {\em geometric genus} is defined to be
$$
g(C)\,:=\,h^0(C,\omega_C)\,=\,h^1(C,\OO_C),
$$
where $\omega_C$ denotes the dualizing sheaf.
The second equality follows from Serre duality.
Since $C$ is smooth over $k$, the sheaf
$\omega_C$ is isomorphic to the sheaf
of K\"ahler differentials $\Omega_{C/k}$.

Let $\varphi: C\to D$ be a finite morphism between smooth curves.
Let us also assume that $\varphi$ is separable, i.e.,
the induced field extension $k(D)\subset k(C)$ is separable.
Then, the {\em Riemann--Hurwitz formula} states that there is
a linear equivalence of divisors on $C$
$$
   K_C \,\lin\, 
   \varphi^*(K_D)\,+\, \sum_{P\in C} {\rm length}(\Omega_{C/D})_P\cdot P\,.
$$
Here, $\Omega_{C/D}$ is the sheaf of relative K\"ahler differentials.
Since $\varphi$ is separable, it is generically \'etale,
and thus, $\Omega_{C/D}$ is a torsion sheaf supported in finitely many points.
By definition, the points in the support of this sheaf are called
{\em ramification points}.
In case $\varphi$ is inseparable, $\Omega_{C/D}$ is non-trivial
in every point, and every point would count as
ramification point.

For a point $P\in C$ with image $Q=\varphi(P)$, and still assuming $\varphi$
to be separable, we choose a local
parameter $t\in\OO_{D,Q}$ and define the {\em ramification index} $e_P$ of
$\varphi$ at $P$ to be the valuation of $\varphi^{\#}(t)$ in $\OO_{C,P}$,
see \cite[Section IV.2]{Hartshorne}. 
Then, the ramification at a ramification point $P$ 
is called {\em tame}, if $e_P$ is not divisible by $p={\rm char}(k)$, 
and it is called {\em wild} otherwise.
We have
$$
  {\rm length}(\Omega_{C/D})_P \,\, \left\{
\begin{array}{cllc}
 = & e_P-1 & \mbox{ \quad if $P$ is tame}\\
 > & e_P-1 & \mbox{ \quad if $P$ is wild}&.
\end{array}
\right.
$$
In general, it is very difficult to bound $e_P$ whenever the ramification
is wild, see the example below.
An important case where one can say more about wild ramification
is in case $\varphi$ is a Galois morphism:
then, one can define for every wild ramification point $P$ certain 
subgroups of the
Galois group, the so-called {\em higher ramification groups},
that control the length of $\Omega_{C/D}$ at $P$, cf.
\cite[Chapitre IV.1]{Serre CorpsLocaux}.

\begin{Example}
  Galois covers with Galois group $\ZZ/p\ZZ$ are called {\em Artin--Schreier
  covers}.
  An example is $\varphi:C\to\PP^1$ given by the 
  projective closure and normalization of the affine equation
  $$
     z^p-z = t^{hp-1}\,.
  $$
  This cover is branched only over $t=\infty$, the ramification is wild 
  of index $e_\infty=p$, and we have ${\rm length}(\Omega_{C/\PP^1})_\infty=p(p-1)h$.
  Thus, $C$ is a curve of genus $1-p+\frac{1}{2}p(p-1)h$ and there are
  $(p-1)h$ non-trivial higher ramification groups.

  In particular, $\varphi$ defines a non-trivial \'etale cover
  of $\Aff^1$, which implies that the affine line $\Aff^1$
  is {\em not} algebraically simply connected.
  In fact, by Raynaud's theorem (formerly Abhyankar's conjecture),
  every finite group that is generated by its $p$-Sylow subgroups
  occurs as quotient of $\piet(\Aff^1)$.
  We refer to \cite{Bost Raynaud} for an overview and references.
  However, it is still true that every \'etale cover
  of $\PP^1$ is trivial, i.e., $\PP^1$ is algebraically
  simply connected \cite[Chapter XI.1]{SGA1}.
\end{Example}

If $\varphi$ is purely inseparable, there still
exists a sort of Riemann--Hurwitz formula.
We refer to \cite{Ekedahl foliations} or \cite[Lecture III]{Miyaoka}
for more information on $\Omega_{C/D}$ in this case.
In this case, the ``ramification divisor'' is defined only up
to linear equivalence.
On the other hand, the structure of purely inseparable morphisms
between curves is simple:
namely, every such morphism is just the composite
of $k$-linear Frobenius morphisms (Proposition \ref{prop:frobenius curves}).
However, from dimension two on, inseparable morphisms become
more complicated.
We will come back to this in Section \ref{sec:Foliations}.
\medskip

Let us now give a rough classification of curves:
If a smooth projective curve over a field $k$
has genus zero, then $\omega_C^\vee$ is
very ample and embeds $C$ as a quadric in $\PP^2_k$.
Moreover, a quadric with a $k$-rational
point is isomorphic to $\PP_k^1$ over any field.
Moreover, the Riemann--Hurwitz formula implies that 
every curve that is dominated by a curve of genus zero
also has genus zero (L\"uroth's theorem).
Thus, since we assumed $k$ to be algebraically closed, we
find

\begin{Theorem}
 If $g(C)=0$ then $C\iso\PP^1_k$, i.e., $C$ is rational. 
 Moreover, every unirational curve, i.e., a curve that is 
 dominated by $\PP^1_k$, is rational.
\end{Theorem}

Although unirational surfaces are rational in characteristic
zero by Castelnuovo's theorem, 
this is false in positive characteristic,
see Theorem \ref{thm:zariski surfaces}.

For curves of genus one, we refer to \cite[Chapter IV.4]{Hartshorne}
or \cite{aec1}. 
Their classification is as follows: 

\begin{Theorem}
 Let $C$ be a smooth projective curve of genus $g(C)=1$ over
 an algebraically closed field $k$ of characteristic $p\geq0$.
 Then
 \begin{enumerate}
   \item after choosing a point $O\in C$ there exists the 
     structure of an Abelian group on the points of $C$,
     i.e., $C$ is an Abelian variety of dimension one,
     an {\em elliptic curve}.
   \item The linear system $|2O|$ defines a finite morphism of degree $2$
     $$  
         \varphi\,:\,C\,\to\,\PP^1_k\,.
     $$
   \item There exists a {\em $j$-invariant} $j(C)\in k$ such that
     two genus one curves are isomorphic if and only if
     they have the same $j$-invariant. 
   \item If $p\neq2$ then $\varphi$ is branched over four points
     and $j$ can be computed from the cross ratio of these
     four points.
   \item The linear system $|3O|$ embeds $C$ as a cubic curve
     into $\PP^2_k$.
     Moreover, if $p\neq2,3$ then $C$
     can be given by an affine equation
     (Weierstra\ss\ equation)
     $$
        y^2\,=\,x^3\,+\,ax\,+\,b
     $$
     for some $a,b\in k$.
 \end{enumerate}
\end{Theorem}

We note that the description of complex elliptic curves
as quotients of $\CC$ by lattices
also has an analog in positive characteristic.
This leads to the theory of {\em Drinfel'd modules},
which is parallel to the theory of elliptic curves
but has not so much to do with
the theory of curves, see \cite[Chapter 4]{Goss}.

A curve $C$ of genus $g\geq2$
is called {\em hyperelliptic} if there exists a 
separable morphism $\varphi:C\to\PP^1_k$ of degree $2$. 
If $p={\rm char}(k)\neq2$ then $\varphi$ is branched over
$2g+2$ points.
(We also note that a morphism of degree $2$
in characteristic $\neq2$ is automatically separable.)
On the other hand, if $p=2$ then every ramification point of $\varphi$
is wildly ramified and thus, there are at most $g+1$ branch
points.
In any characteristic, curves of genus $g=2$ are hyperelliptic
and the generic curve of genus $g\geq3$ is not hyperelliptic.

\begin{Theorem}
 If $g\geq2$ then $\omega_C^{\otimes2}$ is very ample 
 if and only if $C$ is not hyperelliptic.
 In any case, $\omega_C^{\otimes 3}$ is very ample for all curves 
 with $g\geq3$ and $\omega_C^{\otimes 4}$ is very ample for all
 curves of genus $g\geq2$.
\end{Theorem}

Thus, curves of genus $g\geq2$ embed into some fixed projective space 
that depends on $g$ only.
This is the first step towards constructing their 
moduli spaces.
More precisely, Deligne and Mumford \cite{deligne mumford}
showed the existence of a Deligne--Mumford
stack, flat and of dimension $3g-3$ over $\Spec\ZZ$ that
parametrizes curves of genus $g$.
Thus, the moduli space of curves in positive characteristics arises by
reducing the one over $\Spec\ZZ$ modulo $p$.
\medskip

Let us finally mention a couple of facts concerning
automorphism groups:
\begin{enumerate}
 \item If $p\neq2,3$ then the automorphism group
   of an elliptic curve, i.e., automorphisms 
   fixing the neutral element $O$, has order $2$, $4$
   or $6$.
 \item However, the elliptic curve with
   $j=0$ has $12$ automorphisms if $p=3$
   and even $24$ automorphisms if $p=2$, see 
   \cite[Theorem III.10.1]{aec1}.
 \item The automorphism group of a curve of genus $g\geq2$
   is finite. 
   However, the Hurwitz bound $84(g-1)$ on its
   in characteristic zero can be violated.
   We refer to \cite[Chapter IV.2, Exercise 2.5]{Hartshorne} 
   for details and further references
\end{enumerate}

Let us note that some classes of surfaces arise as 
quotients $(C_1\times C_2)/G$, where $C_1$, $C_2$ are curves with
$G$-actions.
Now, in positive characteristic larger automorphism groups 
may show up and thus, new possibilities have to be considered.
For example, we will see in Section \ref{subsec:hyperelliptic surfaces} that
hyperelliptic surfaces arise as quotients of products of
elliptic curves in any characteristic.
It is remarkable that {\em no} new classes arise
in characteristic $2$ and $3$ from larger
automorphism groups of elliptic curves with $j=0$.

\subsection{Group schemes}
\label{subsec:group schemes}
Constructions with groups are ubiquitous in geometry.
Instead of finite groups we will consider {\em finite and flat group schemes}
$G$ over a ground field $k$, which we assume to be algebraically
closed of characteristic $p\geq0$.
We refer to \cite{Waterhouse} or \cite{Tate group schemes}
for overview, details and references.

Thus, $G=\Spec A$ for some finite-dimensional $k$-algebra $A$,
and there exist morphisms
$$
O\,:\,\Spec k\to G \mbox{ \quad and \quad }
m\,:\,G\times G\to G
$$
where $m$ stands for multiplication, and $O$ for the neutral element.
These morphisms have to fulfill certain axioms that encode that $G$ 
is a group object in the category of schemes.
We refer to \cite[Chapter I]{Waterhouse} for the precise definition and
note that it amounts to saying that $A$ is
a {\em Hopf algebra}.
The dimension $\dim_k A$ is called the {\em length}, or {\em order}, 
of the group scheme $G$.

The following construction associates to every finite group
a finite flat group scheme:
for a finite group or order $n$ with elements $g_1,...,g_n$
we take a disjoint union of
$n$ copies of $\Spec k$, one representing each $g_i$,
and define $m$ via the multiplication in the group we started with.
This defines the {\em constant group scheme} associated to a finite group.
Conversely, we have 

\begin{Theorem}
 A finite flat group scheme $G$ 
 of length prime to $p$ over an algebraically closed field
 is a constant group scheme.
\end{Theorem}

In particular, over an algebraically closed field of characteristic zero,
we obtain an equivalence
between the categories of finite groups and finite flat group schemes.

One feature of constant group schemes is that the structure morphism
$G\to\Spec k$ is \'etale, i.e., $A=H^0(G,\OO_G)$ is a separable $k$-algebra.
For example, consider the constant group scheme $\ZZ/p\ZZ$, which is of length $p$.
As an algebra, $A$ is isomorphic to $k^p$ with componentwise
addition and multiplication, and thus \'etale over $k$.
On the other hand, we will see below that there are two different
structures of group schemes on $\Spec k[x]/(x^p)$, which is not reduced -- these are 
examples of {\em infinitesimal group schemes}.
Let us first note that infinitesimal group schemes are a particular 
characteristic $p$ phenomenon:

\begin{Theorem}[Cartier]
 \label{thm:cartier}
 Group schemes over fields of characteristic zero are smooth and thus, reduced.
\end{Theorem}

To give examples of infinitesimal group schemes, we consider
$\GG_a$ and $\GG_m$.
Here, $\GG_a$ denotes the group scheme corresponding to the additive group,
i.e., $(\GG_a(k),\circ)=(k,+)$.
Similarly, $\GG_m$ denotes the group scheme corresponding to the 
multiplicative group of $k$, 
i.e., $(\GG_m(k),\circ)=(k^\times,\cdot)$,
see \cite[Chapter I.2]{Waterhouse} -
these group schemes are affine but not finite over $k$.
Then, the first example of an infinitesimal group scheme
is $\mu_p$, the group scheme of {\em $p$.th roots of unity}.
Namely, there exists a short exact sequence of group schemes
(in the flat topology)
$$
0\,\to\,\mu_p\,\to\,\GG_m\,\stackrel{x\mapsto x^p}{\longrightarrow}\,
\GG_m\,\to\,0\,.
$$
We note that the kernel $\mu_p$ is infinitesimal because of the equality 
$x^p-1=(x-1)^p$ in characteristic $p$.
The second example is $\alpha_p$, the kernel of Frobenius on $\GG_a$, i.e.,
we have a short exact sequence
$$
0\,\to\,\alpha_p\,\to\,\GG_a\,\stackrel{F}{\longrightarrow}\,
\GG_a\,\to\,0\,.
$$
Both group schemes, $\alpha_p$ and $\mu_p$, are isomorphic to
$\Spec k[x]/(x^p)$ as schemes, and are thus infinitesimal (non-reduced),
but have different multiplication maps.
Together with $\ZZ/p\ZZ$ these are all group schemes of length $p$:

\begin{Theorem}[Oort--Tate]
\label{thm:oort tate}
  A finite flat group scheme of length $p$ over an  
 algebraically closed field of characteristic $p$ is isomorphic to
 $\ZZ/p\ZZ$, $\mu_p$ or $\alpha_p$.
\end{Theorem}

For more general results we refer to \cite{Tate Oort}
and \cite{oort}.
Let us also mention that there exist non-Abelian group schemes of order
$p^2$.  
Thus, in positive characteristic, the theory of finite flat group schemes is richer than
the theory of finite groups, already over an algebraically closed field.

For example, if $E$ is an elliptic curve in characteristic
$p$, then multiplication by $p$ induces a morphism $E\to E$,
whose kernel $E[p]$ is a finite flat group scheme of length
$p^2$ (as expected from characteristic zero).
More precisely, and still assuming $k$ to be algebraically
closed,
$$
E[p] \,\iso\,\left\{ 
\begin{array}{lll}
 \mbox{ either } &
 \mu_p\oplus(\ZZ/p\ZZ) & \mbox{ and $E$ is called {\em ordinary},}\\
 \mbox{ or } &
 M_2 & \mbox{ a non-split extension of $\alpha_p$ by itself}\\
 & & \mbox{ and $E$ is called {\em supersingular}.}
\end{array}
\right.
$$
Looking at $k$-rational points, we find 
$E[p](k)=\ZZ/p\ZZ$ if $E$ is ordinary, and
$E[p](k)=0$ if $E$ is supersingular.
Thus, $k$-rational points do not suffice to see the full $p$-torsion,
and the theory of finite flat group schemes is really needed.
As the name suggests, the generic elliptic curve is ordinary.
More precisely, a theorem of Deuring states that
there exist approximately $p/12$ supersingular elliptic curves 
in characteristic $p$, see \cite[Theorem V.4.1]{aec1}.

In classical algebraic geometry,
one often constructs interesting and new varieties as Galois-covers or quotients by finite
groups of ``well-understood'' varieties.
In positive characteristic, one very successful way to construct a 
``pathological characteristic-$p$''
example is via purely
inseparable covers, or, via quotients by
infinitesimal group schemes.
The role of Galois covers is often played by torsors under
$\alpha_p$ and $\mu_p$.
We come back to this in Section \ref{subsec:quot by alpha}.

\section{Cohomological tools and invariants}
\label{sec:Cohomology}

This section circles around algebraic 
versions of Betti and Hodge numbers, and deRham-cohomology.
Especially towards the end, the subjects get deeper, our
exposition becomes sketchier and we
advise the reader interested in surface theory only, to skip
all but the first three paragraphs.

In this section, $X$ will be a smooth and projective
variety of arbitrary dimension over an algebraically
closed field $k$ of characteristic $p\geq0$.

\subsection{Hodge numbers}
As usual, we define the {\em Hodge numbers} of $X$ to be
$$
h^{i,j}(X)\,:=\,\dim_k H^j(X,\Omega_X^i)\,.
$$
We note that Serre duality holds for projective 
Cohen--Macaulay schemes over any field
\cite[Chapter III.7]{Hartshorne}, and in 
particular we find
$$
   h^{i,j}(X) \,=\, h^{n-i, n-j}(X),
   \mbox{ \quad where \quad }n=\dim(X)\,.
$$
Over the complex numbers, complex conjugation induces the
Hodge symmetry $h^{i,j}=h^{j,i}$, see, for example 
\cite[Chapter 0.7]{Griffiths Harris}.
However, even for a smooth projective surface in 
positive characteristic, the
numbers $h^{0,1}$ and $h^{1,0}$ may be 
different.
For example, in \cite[Theorem 8.3]{Liedtke Uniruled}
we constructed a sequence $\{X_i\}_{i\in\NN}$ of surfaces 
with fixed $\Pic^0_{\rm red}$
in characteristic $2$, where 
$h^{1,0}(X_i)-h^{0,1}(X_i)$ 
tends to infinity.

\subsection{Betti numbers}
An algebraic replacement for singular cohomology is
{\em $\ell$-adic cohomology}, whose construction is
due to Grothendieck.
We refer to \cite[Appendix C]{Hartshorne} for motivation,
as well as to \cite{Milne} for a complete
treatment.
Let us here only describe its basic properties:
let $\ell$ be a prime number different from $p$
and let $\QQ_\ell$ be the field of
$\ell$-adic numbers, i.e., the completion of
$\QQ$ with respect to the $\ell$-adic valuation,
see also Section \ref{subsec:Witt def}.
Then, 
\begin{enumerate}
 \item the $\ell$-adic cohomology groups 
   $\Het{i}(X,\QQ_l)$ are finite-dimensional
   $\QQ_\ell$-vector spaces,
 \item they are zero for $i<0$ and $i>2\dim(X)$,
 \item the dimension of $\Het{i}(X,\QQ_\ell)$
  is independent of $\ell$ (here, 
  $\ell\neq p$ is crucial), and we denote it by
  $b_i(X)$, the $i$.th {\em Betti number},
 \item $\Het{*}(X,\QQ_\ell)$ satisfies Poincar\'e duality.
\end{enumerate}

If $k=\CC$, then so-called {\em comparison theorems} show 
that these Betti numbers coincide with the topological ones.
Let us also mention the following feature:
if $k$ is not algebraically closed, then the absolute
Galois group ${\rm Gal}(\overline{k}/k)$ acts on the
$\ell$-adic cohomology groups of $X_{\overline{k}}$, 
which gives rise to interesting representations of ${\rm Gal}(\overline{k}/k)$.

For the following two classes of varieties,
$\ell$-adic cohomology and Hodge invariants
are precisely as one would expect them from complex geometry:
\begin{enumerate}
  \item If $C$ is a smooth and projective curve over $k$, then
   $b_0=b_2=1$, $b_1=2g$, and $h^{1,0}=h^{0,1}=g$. 
  \item If $A$ is an Abelian variety of dimension $g$ over 
   $k$ then $b_0=b_{2g}=1$, $b_1=2g$, and
   $h^{0,1}=h^{1,0}=g$.
   Moreover, there exists an isomorphism
   $$
    \Lambda^i \Het{1}(A,\QQ_l)\,\iso\,\Het{i}(A,\QQ_\ell)
   $$
   giving -- among many other things -- the expected Betti numbers.
\end{enumerate}

However, for more general classes of smooth and projective
varieties,
the relations between Betti numbers, Hodge invariants, 
deRham-cohomology and the Fr\"olicher spectral 
sequence in positive characteristic
are more subtle than over the complex numbers,
as we shall see below.

Let us first discuss $h^{1,0}$, $h^{0,1}$ and $b_1$ 
in more detail, since this is important for the classification of
surfaces.
Also, these numbers can be treated fairly elementary.

\subsection{Picard scheme and Albanese variety}
If $X$ is smooth and proper over a field $k$, then there
exists an Abelian variety ${\rm Alb(X)}$ over $k$, the 
{\em Albanese variety} of $X$, and an
{\em Albanese morphism}
$$
{\rm alb}_X\,:\,X\,\to\,\Alb(X)\,.
$$
The pair $(\Alb(X),{\rm alb}_X)$ is characterized by 
the universal property that every morphism from $X$ to
an Abelian variety factors over ${\rm alb}_X$.
For a purely algebraic construction, we refer to 
\cite{Serre Albanese}.

Next, the Picard functor, which classifies invertible sheaves on $X$, is
representable by a group scheme, the {\em Picard scheme} $\Pic(X)$,
whose neutral element is $[\OO_X]$, see \cite{FGA Picard}.
We denote by $\Pic^0(X)$ the identity component of $\Pic(X)$.
Deformation theory provides us with a natural
isomorphism 
$$
T\Pic^0(X)\,\iso\,H^1(X,\OO_X)\,,
$$
where $T\Pic^0(X)$ denotes the Zariski-tangent space at 
$[\OO_X]$.

Now, group schemes over fields of
positive characteristic may be non-reduced 
(the group schemes $\mu_p$ and $\alpha_p$ from 
Section \ref{subsec:group schemes} are examples), but 
the reduction of $\Pic^0(X)$ is still an Abelian variety, 
which is the dual Abelian variety of $\Alb(X)$,
see \cite[Chapter 5]{Badescu}.
Also, the first Betti number $b_1$ is twice the dimension
of $\Alb(X)$. 
Thus, we get
$$
\frac{1}{2}b_1(X)\,=\,\dim \Alb(X)\,=\,\dim\Pic^0(X).
$$
Since, the dimension of the Zariski tangent space
at $[\OO_X]\in\Pic^0(X)$ is at least equal to the dimension 
of $\Pic^0(X)$, we find 
$$
h^{0,1}(X)\,=\,h^1(X,\OO_X)\,\geq\,\frac{1}{2}b_1(X)\,,
$$
with equality if and only if $\Pic^0(X)$ is a reduced group scheme, 
i.e., an Abelian variety.
By Cartier's theorem (Theorem \ref{thm:cartier}), group schemes over a 
field of characteristic zero are reduced.
As a corollary, we obtain a purely algebraic proof of the
following fact 

\begin{Proposition}
 A smooth and proper variety over a field of characteristic zero
 satisfies $b_1(X)=2h^{0,1}(X)$.
\end{Proposition}

For curves and Abelian varieties over arbitrary fields,
$b_1$, $h^{1,0}$ and $h^{0,1}$ are precisely as
over the complex numbers.
On the other hand, over fields of positive characteristic,
\begin{enumerate}
\item  there do exist surfaces with $h^{0,1}>b_1/2$, i.e., 
  with non-reduced  Picard schemes, 
  see \cite{Igusa OpenProblems} and \cite{Serre Mexico}.
\end{enumerate}

In \cite[Lecture 27]{Mumford Curves},
the non-reducedness of $\Pic^0(X)$ is related
to non-trivial Bockstein operations 
$\beta_n$ from subspaces of $H^1(X,\OO_X)$ to
quotients of $H^2(X,\OO_X)$.
In particular, a smooth projective variety with
$h^2(X,\OO_X)=0$ has a reduced $\Pic^0(X)$,
which applies, for example, to curves.
In the case of surfaces, a quantitative 
analysis of which classes can
have non-reduced Picard schemes
has been carried out in \cite{Liedtke Picard}.

\subsection{Differential one-forms}
\label{subsec:one-forms}
We shall see in Section \ref{subsec:instructive example} that in
positive characteristic, the
pull-back of a non-zero differential form under a morphism
may become zero.
However, by a fundamental theorem of Igusa \cite{Igusa Fundamental}, 
every non-trivial global $1$-form on $\Alb(X)$ pulls back, 
via ${\rm alb}_X$, to a non-zero global $1$-form on $X$.
This implies the estimate
$$
  h^{1,0}(X)\,=\, h^0(X,\Omega_X^1)\,\geq\,\frac{1}{2}b_1(X)\,.
$$
Moreover, all global $1$-forms arising as pull-back from $\Alb(X)$
are $d$-closed, i.e., closed under the exterior derivative.
Despite of Igusa's theorem, it is still possible that the Albanese morphism in positive characteristic 
becomes purely inseparable (see Theorem \ref{thm:purely insep Albanese}).

We have $h^{1,0}=b_1/2$ for curves and Abelian varieties 
over arbitrary fields and their  global $1$-forms are $d$-closed.
On the other hand, over fields of positive characteristic,
\begin{enumerate}
 \item there do exist surfaces 
   with $h^{1,0}>b_1/2$, i.e., with 
   ``too many'' global $1$-forms, see 
   \cite{Igusa OpenProblems}, and
 \item there do exist surfaces with global $1$-forms that are not
   $d$-closed, see \cite{Mumford Path1} and \cite{Fossum}. 
   These forms give rise to a non-zero differential in their
   Fr\"olicher spectral sequences, which thus do {\em not}
   degenerate at $E_1$.
\end{enumerate}
We refer to \cite[Proposition II.5.16]{Illusie deRham-Witt} 
for more results and to \cite[Section II.6.9]{Illusie deRham-Witt} 
for the connection to Oda's subspace in first deRham 
cohomology.

\subsection{Igusa's inequality}
\label{subsec:Igusa}
We denote  by $\rho$ the rank of the N\'eron--Severi group ${\rm NS}(X)$,
which is always finite.
More precisely, Igusa's theorem \cite{Igusa Betti} states
$$
\rho(X)\,\leq\,b_2(X)\,.
$$
This follows from the existence of a Chern map from $\NS(X)$ to
second $\ell$-adic or crystalline cohomology.
On the other hand, $d\log$ induces a ``naive'' cycle map
$$
d\log\,:\,\NS(X)\otimes_\ZZ k \,\to\, H^1(X,\Omega_X^1)\,,
$$
which is injective in characteristic zero, and which then implies 
the inequality $\rho\leq h^{1,1}$,
see, for example, \cite[Exercise 5.5]{Badescu}.
However, this map may fail to be injective in positive characteristic, as
the example of supersingular Fermat surfaces 
\cite{Shioda 1974} shows.
More precisely, these surfaces 
satisfy $b_2=\rho>h^{1,1}$, see also Section \ref{subsec:fermat}.

\subsection{Kodaira vanishing}
\label{subsec:kodaira fails}
Raynaud \cite{ray} gave the first counter-examples  
to the Kodaira vanishing theorem in
positive characteristic.
However, we mention the following results that
tell us that the situation is not too bad:
\begin{enumerate}
 \item If ${\cal L}$ is an ample line bundle then
  ${\cal L}^{\otimes \nu}$, $\nu\gg0$ fulfills Kodaira
  vanishing (in fact, this is just Serre vanishing) 
  \cite[Theorem III.7.6]{Hartshorne}. 
 \item  If a smooth projective variety 
  of dimension $<p$ lifts over $W_2(k)$ then Kodaira vanishing holds, 
  see \cite[Theorem 5.8]{Illusie Frobenius},
  \cite[Corollaire 2.11]{Deligne Illusie}, and Section \ref{subsec:witt lift}.
  Under stronger lifting assumptions, also Kawamata--Viehweg
  vanishing holds \cite{Xie vanishing}.
 \item Kodaira vanishing, and even stronger vanishing results hold
  for the (admittedly rather special) class of 
  {\em Frobenius-split} varieties \cite[Theorem 1.2.9]{Brion Kumar}.
 \item In \cite[Section II]{ek}, Ekedahl develops tools to handle 
  possible failures of Kodaira vanishing, see also
  Section \ref{subsec:pluricanonical}.
 \item In \cite{Xie converse}, Xie shows that all surfaces violating
  Kodaira-Ramanujam vanishing arise as in \cite{ray}.
\end{enumerate}
Here, $W(k)$ denotes the ring of Witt vectors and 
$W_2(k)$ the ring of Witt vectors of length $2$,
see Section \ref{subsec:Witt def}.
Let us just mention that if $k$ is a perfect field 
then $W(k)$ is a complete discrete valuation ring of
characteristic zero with residue field $k$, and that
this ring is in a certain sense minimal and universal.

\subsection{Fr\"olicher spectral sequence}
Let $\Omega_X^i$ be the sheaf of (algebraic) 
differential $i$-forms.
These sheaves, together with the exterior derivative $d$
form a complex, the
{\em (algebraic) deRham-complex} $(\Omega_X^*,d)$.
Now, the Zariski topology is too coarse to have a
Poincar\'e lemma.
Thus, we define {\em (algebraic) deRham-cohomology} $\HdR{*}(X/k)$ 
to be the hypercohomology of this complex.
In particular, there always exists a spectral sequence
$$
E_1^{i,j}\,=\,H^j(X,\Omega_X^i) \,\Rightarrow\, \HdR{i+j}(X/k),
$$
the {\em Fr\"olicher spectral sequence}, from Hodge- to
deRham-cohomology.
If $k=\CC$, and $X$ is proper over $k$,
then these cohomology groups and the spectral
sequence coincide with the analytic ones, see
\cite{Grothendieck deRham}.
Already the existence of the Fr\"olicher spectral sequence 
implies for all $m\geq0$ the inequality
$$
    \sum_{i+j=m} h^j(X,\Omega^i_X)\,\geq\, \hdR{m}(X/k)\,.
$$
Equality for all $m$ is equivalent to the degeneration of this 
spectral sequence at $E_1$.
Over the complex numbers, degeneration at $E_1$ is true - 
however, the classical proof uses methods from differential geometry,
functional analysis and partial differential equations,
see \cite[Chapter 0.7]{Griffiths Harris}.
On the other hand, if a variety of positive characteristic
admits a lift over $W_2(k)$, 
then we have the following result from
\cite{Deligne Illusie} (but see
\cite{Illusie Frobenius} for an expanded version):

\begin{Theorem}[Deligne--Illusie]
 Let $X$ be a smooth and projective variety 
 in characteristic $p\geq\dim X$ and assume
 that $X$ admits a lift over $W_2(k)$.
 Then the Fr\"olicher spectral sequence of $X$
 degenerates at $E_1$.
\end{Theorem}

The assumptions are fulfilled for curves and Abelian
varieties, see Section \ref{sec:Witt}.
Moreover, if a smooth projective variety $X$ 
in characteristic zero admits a model over $W(k)$ 
for some perfect field of characteristic 
$p\geq\dim X$ it follows from semi-continuity
that the Fr\"olicher spectral sequence of $X$
degenerates at $E_1$ in characteristic zero.
From this one obtains purely algebraic proofs of 
the following

\begin{Theorem}
  Degeneration 
  at $E_1$ holds for
 \begin{enumerate}
  \item smooth projective curves and Abelian varieties over arbitrary fields, and
  \item smooth projective varieties over fields of 
    characteristic zero. 
 \end{enumerate}
\end{Theorem}

We already mentioned above that varieties with 
global $1$-forms that are not $d$-closed, such as Mumford's surfaces 
\cite{Mumford Path1}, provide examples where
degeneration at $E_1$ does {\em not} hold.

\subsection{Crystalline cohomology}
\label{subsec:crystalline}
To link deRham-, Betti- and
Hodge-cohomology, we use {\em crystalline cohomology}.
Its construction, due to Berthelot and Grothendieck,
is quite involved \cite{Berthelot}.
This cohomology theory takes values in the Witt
ring $W=W(k)$, which is a discrete valuation ring if
$k$ is perfect, see Section \ref{sec:Witt}.
In case a smooth projective variety lifts to 
some ${\cal X}/W(k)$,
crystalline cohomology is the deRham cohomology
$\HdR{*}({\cal X}/W(k))$.
It was Grothendieck's insight, and the starting point
of crystalline cohomology, that this deRham cohomology
does not depend on the choice of lift ${\cal X}$.
One of the main technical difficulties to overcome defining
crystalline cohomology for arbitrary smooth and proper
varieties is that they usually do {\em not} lift over
$W(k)$.

If $X$ is a smooth projective variety over a perfect field $k$
then
\begin{enumerate}
 \item the groups $\Hcris{i}(X/W)$ are finitely
   generated $W$-modules,
 \item they are zero for $i<0$ and $i>2\dim(X)$,
 \item there are actions of Frobenius and Verschiebung on $\Hcris{i}(X/W)$,
 \item $\Hcris{*}(X/W)\otimes_W K$ satisfies Poincar\'e duality, where
  $K$ denotes the field of fractions of $W$, and
 \item if $X$ lifts over $W(k)$ then crystalline cohomology is
  isomorphic to the deRham cohomology of a lift.
\end{enumerate}
We remind the reader that in order to get the ``right'' Betti numbers
from the $\ell$-adic cohomology groups, 
we had to assume $\ell\neq p$.
Crystalline cohomology takes values in $W(k)$ 
(recall $W(\FF_p)\iso\ZZ_p$ with field of fractions $\QQ_p$), 
and this is the ``right'' cohomology theory for $\ell=p$. 
In fact,
$$
   b_i(X)\,=\,\dim_{\QQ_\ell} \Het{i}(X,\QQ_\ell) 
   \,=\, {\rm rank}_W\, \Hcris{i}(X/W)\mbox{ \quad for all \quad }\ell\neq p,
$$
i.e., the Betti numbers of $X$ are encoded in the rank of
crystalline cohomology.
However, since the $\Hcris{*}(X/W)$ are $W$-modules, there may be
non-trivial torsion - and this is precisely the explanation
for the differences between Hodge- and Betti-numbers.
More precisely, there is a universal coefficient formula,
and for all $m\geq0$ there are short exact sequences
$$
0\,\to\,\Hcris{m}(X/W)\otimes_W k\,\to\,\HdR{m}(X/k)\,\to\,
{\rm Tor}^{W}_1(\Hcris{m+1}(X/W),k)\,\to\,0\,.
$$
(In view of what we already said in case $X$ admits a lift over $W$,
it should be plausible that there is a connection between crystalline 
and deRham cohomology.)
In particular, Betti- and deRham-numbers coincide if and only
if all crystalline cohomology groups are torsion-free
$W$-modules.

\subsection{Hodge--Witt cohomology}
In \cite{Illusie deRham-Witt}, Illusie constructed 
the {\em deRham-Witt complex} $W\Omega_X^*$ and studied its
cohomology groups $H^j(X,W\Omega_X^i)$, the 
{\em Hodge--Witt cohomology groups}.
For $i=0$, these coincide with Serre's Witt vector cohomology groups 
introduced in \cite{Serre Mexico}.
The Hodge--Witt cohomology groups are
$W$-modules, whose torsion may not be finitely generated.
In any case, there exists a spectral sequence, the 
{\em slope spectral sequence}
$$
   E_1^{i,j}\,=\,H^j(X,W\Omega_X^i)\,\Rightarrow\,\Hcris{i+j}(X/W)\,,
$$
which degenerates at $E_1$ modulo torsion.
We refer to \cite[Section II.7]{Illusie deRham-Witt} for computations and 
further results.

Finally, using Hodge--Witt cohomology and slopes on crystalline cohomology,
Ekedahl \cite[page 85]{Ekedahl Gauge} (but see also \cite{Illusie Ekedahl})
proposed new invariants of smooth
projective varieties: slope numbers, dominoes and Hodge--Witt numbers.
It is not yet clear, what role they will eventually play in characteristic-$p$ geometry. 
We refer to \cite{Joshi} for some results.

\addtocontents{toc}{{\medskip \bf Classification of Algebraic Surfaces}}

\section{Birational geometry of surfaces}
\label{sec:Birational}

From this section on, we study smooth surfaces.
To start with, we discuss their birational geometry,
which turns out to be ``basically the same'' as over the
complex numbers.
Unless otherwise stated, results and proofs can be found in
\cite{Badescu}, and we refer to \cite{BoHu} for an
overview different from ours.

\subsection{Riemann--Roch}
\label{subsec:riemann-roch}
Let $S$ be a smooth projective surface over an
algebraically closed field $k$ of characteristic $p\geq0$.
Actually, asking for properness would be enough:
by a theorem of Zariski and Goodman, a surface that is
smooth and proper over a field is automatically projective,
see \cite[Theorem 1.28]{Badescu}.

For every locally free sheaf $\cal E$,
Grothendieck constructed Chern-classes $c_i({\cal E})$ that take
values in Chow-groups, $\ell$-adic, or crystalline cohomology.
As usual, for a smooth variety $X$ with tangent sheaf $\Theta_X$
we set $c_i(X):=c_i(\Theta_X)$.

We have {\em Noether's formula}
$$
\chi(\OO_S) \,=\, \frac{1}{12}\left(\, c_1^2(S)\,+\,c_2(S)\,\right)\,.
$$
Moreover, if ${\cal L}$ is an invertible sheaf on $S$,
we have the {\em Riemann--Roch formula} 
$$
\chi({\cal L}) \,=\, \chi(\OO_S)\,+\,\frac{1}{2}{\cal L}\cdot ({\cal L}\otimes\omega_S^\vee)\,.
$$
We note that Serre duality holds 
for Cohen--Macaulay schemes that are of finite type over a field.
Thus, we have an equality
$h^i(S,{\cal L})=h^{2-i}(S,\omega_S\otimes{\cal L}^\vee)$ for
surfaces.
However, we have seen in Section \ref{subsec:kodaira fails} that 
Kodaira vanishing may not hold.
Finally, if $D$ is an effective divisor on $S$, then $D$ is
a Gorenstein curve and the {\em adjunction formula} yields
$$
\omega_D \,\iso\, \left(\omega_S\otimes\OO_S(D)\right)|_D\,,
$$
where $\omega_D$ and $\omega_S$ denote the respective dualizing sheaves.
In particular, if $D$ is reduced and irreducible, we obtain
$$
2p_a(D)\,-\,2\,=\,D^2\,+\,K_S \cdot D\,,
$$
where $p_a$ denotes the arithmetic genus of $D$.
We refer to \cite[Chapter V.1]{Hartshorne}, 
\cite[Appendix A]{Hartshorne}, \cite[Chapter 5]{Badescu}, 
\cite{Fulton}, and \cite{Milne} for details 
and further references.

\subsection{Blowing up and down}
First of all, blowing up a point on a smooth surface over an
algebraically closed field has the same effect as over the complex numbers.

\begin{Proposition}
 Let $f:\widetilde{S}\to S$ be the blow-up in a closed point 
 and denote by $E$ the exceptional divisor on $\widetilde{S}$.
 Then
 $$
    E\,\iso\,\PP^1_k,\,\mbox{ \quad }E^2=-1,\mbox{ \quad and \quad }K_{\widetilde{S}} \cdot E=-1\,.
 $$
 Moreover, the equalities
 $$
   b_2(\widetilde{S}) \,=\, b_2(S)\,+\,1 \mbox{ \quad and \quad }
   \rho(\widetilde{S}) \,=\,\rho(S)\,+\,1
 $$
 hold true.
\end{Proposition}

As in the complex case, we call such a curve $E$ with
$E^2=-1$ and $E\iso\PP^1$ an {\em exceptional $(-1)$-curve}.
A surface that does not contain exceptional $(-1)$-curves
is called {\em minimal}.

Conversely, exceptional $(-1)$-curves can be contracted and 
the proof (modifying a suitable hyperplane section)
is basically the same as in characteristic zero, cf.
\cite[Theorem 3.30]{Badescu} or
\cite[Theorem V.5.7]{Hartshorne}

\begin{Theorem}[Castelnuovo]
  Let $E$ be an exceptional $(-1)$-curve on a smooth surface
  $S$.
  Then, there exists a smooth surface $S'$ and
  a morphism $f:S\to S'$, such that $f$ is the blow-up
  of $S'$ in a closed point with exceptional divisor $E$.
\end{Theorem}

Since $b_2$ drops every time one contracts an exceptional
$(-1)$-curve, Castelnuovo's theorem implies
that for every surface $S$ there exists a sequence of
blow-downs $S\to S'$ onto a minimal surface $S'$.
In this case, $S'$ is called a {\em minimal model}
of $S$.

\subsection{Resolution of indeterminacy}
As in characteristic zero, a rational map from a surface extends to
a morphism after a finite number of blow-ups in closed points, which
gives {\em resolution of indeterminacy} of a rational map.
Moreover, every birational (rational) map can be factored as
a sequence of blow-ups and Castelnuovo blow-downs, see, e.g.,
\cite[Chapter V.5]{Hartshorne}.

\subsection{Kodaira dimension}
As over the complex numbers, the following two notions are crucial
in the Kodaira--Enriques classification of surfaces:
first, the $n$.th {\em plurigenus} $P_n(X)$ of a smooth projective variety
$X$ is defined to be
$$
    P_n(X) \,:=\, h^0(X,\omega_X^{\otimes n}).
$$
Second, the {\em Kodaira dimension} $\kappa(X)$ is defined to be
$\kappa(X)=-\infty$ if $P_n(X)=0$ for all $n\geq1$, or else
$$ 
   \kappa(X)\,:=\,\max_{n\in\NN}\left\{\, \dim\phi_n(X)\, \right\}
$$
where $\phi_n:X\dashrightarrow \PP^{P_n(X)-1}_k$ denotes
the $n$.th pluri-canonical (possibly only rational) map.
This recalled, we have the following important result,
cf. \cite[Corollary 10.22]{Badescu}, which is already 
non-trivial in characteristic zero:

\begin{Theorem}
  Let $S$ be a smooth projective surface with $\kappa(S)\geq0$.
  Then, $S$ possesses a unique minimal model.
\end{Theorem}

\subsection{Birationally ruled surfaces}
We recall that a surface $S$ is called {\em birationally ruled}
if it is birational to 
$\PP^1\times C$ for some smooth curve $C$.
Such surfaces are easily seen to
satisfy $P_n(S)=0$ for all $n\geq1$, i.e.,
they are of Kodaira dimension $\kappa(S)=-\infty$.
Conversely, one can show (see, e.g., \cite[Theorem 13.2]{Badescu})
that such surfaces with $\kappa(S)=-\infty$ 
possess a smooth rational curve 
that moves:

\begin{Theorem}
 \label{birationally ruled}
 If $S$ is birationally ruled, then $\kappa(S)=-\infty$.
 Conversely, if $\kappa(S)=-\infty$ then $S$ is birationally ruled, i.e.,
 birational to $\PP^1\times C$, and
 $$
   h^1(S,\OO_S)\,=\,\frac{1}{2}b_1(S)\,=\,g(C)\,,
 $$
 where $g(C)$ denotes the genus of $C$.
\end{Theorem}

As in the complex case, minimal models for surfaces with
$\kappa(S)=-\infty$ are not unique.
More precisely, we have Nagata's result

\begin{Theorem}
 \label{thm:nagata}
  Let $S$ be a minimal surface with $\kappa(S)=-\infty$.
  \begin{enumerate}
   \item If $h^1(S,\OO_S)\geq1$, then the image $C$ of the
     Albanese map is a smooth curve.
     Moreover, there exists a rank two vector bundle ${\cal E}$ on
     $C$ such that  ${\rm alb}_S:S\to C$ is isomorphic 
     to $\PP({\cal E})\to C$.
   \item If $h^1(S,\OO_S)=0$, then $S$ is isomorphic to $\PP^2$
      or a Hirzebruch surface 
      $\FF_d:=\PP(\OO_{\PP^1}\oplus\OO_{\PP^1}(d))\to\PP^1$
      with $d\neq1$. 
  \end{enumerate}
\end{Theorem}

Also, Castelnuovo's cohomological characterization of rational
surfaces holds true. 
The proof in positive characteristic
is due to Zariski \cite{Zariski}, but see also the discussion
in \cite[Part 4]{BoHu}:

\begin{Theorem}[Castelnuovo--Zariski]
  \label{thm:castelnuovo zariski}
  For a smooth projective surface $S$, the following are equivalent
  \begin{enumerate}
   \item $S$ is rational, i.e., birational to $\PP^2$,
   \item $h^1(S,\OO_S)=P_2(S)=0$
   \item $b_1(S)=P_2(S)=0$.
  \end{enumerate}
\end{Theorem}

So far, things look pretty much the same as over the
complex numbers.
However, one has to be a little bit careful with the notion
of uniruledness:
we will see in Section \ref{uni and super} below that
unirationality (resp. uniruledness) does {\em not} imply rationality
(resp. ruledness).

\subsection{Del~Pezzo surfaces}
A surface $S$ is called {\em del~Pezzo}, or, {\em Fano},
if $\omega_S^\vee$ is ample.
In every characteristic, these surfaces are rational.
More precisely, they are isomorphic to $\PP^1\times\PP^1$,
$\PP^2$, or $\PP^2$ blown-up in at most $8$ points in general position.
We refer to V\'arilly-Alvarado's lecture notes in this volume for
details.

\section{(Quasi-)elliptic fibrations}
\label{sec:Elliptic}

For the classification of surfaces of special type
in characteristic zero, elliptic fibrations play an
important role.
In positive characteristic, {\em wild fibers}
and {\em quasi-elliptic fibrations} are new features
that show up.
We refer to \cite[Section 7]{Badescu} for an introduction,
\cite{bm2}, \cite{bm3} for more details, to 
\cite{Schuett Shioda} for an overview article with many examples,
and to
\cite[Chapter V]{Cossec; Dolgachev} for more advanced topics.

\subsection{Quasi-elliptic fibrations}
Given a dominant morphism from a smooth 
surface $S$ onto
a curve in any characteristic, we may pass to its 
Stein factorization and obtain a fibration
$S\to B$, cf. \cite[Corollary III.11.5]{Hartshorne}.
Then, its generic fiber $S_\eta$ is an integral curve, i.e.,
reduced and irreducible \cite[Theorem 7.1]{Badescu}.
Moreover, in characteristic zero
Bertini's theorem implies that $S_\eta$
is in fact smooth over the function field
$k(B)$.
Now, if ${\rm char}(k)=p>0$, then it is still true that
the generic fiber is a regular curve, i.e., all
local rings are regular local rings.
However, this does {\em not} necessarily imply that 
$S_\eta$ is smooth over $k(B)$.
Note that since $S_\eta$ is one-dimensional,
regularity is the same as normality.
We refer to \cite[Chapter 11.28]{Matsumura} for a 
discussion of smoothness versus regularity.

Suppose $S_\eta$ is not smooth over $k(B)$ and denote
by $\overline{k(B)}$ the algebraic closure
of $k(B)$.
Then $S_{\overline{\eta}}:=S_\eta\times_{\Spec k(B)}\Spec\overline{k(B)}$
is still reduced and irreducible \cite[Theorem 7.1]{Badescu}
but no longer regular and we denote by
$S_{\overline{\eta}}'\to S_{\overline{\eta}}$ its
normalization.
Then, Tate's theorem on genus change in inseparable extensions
\cite{Tate} 
(see \cite{Schroeer 2009} for a modern treatment)
states

\begin{Theorem}[Tate]
 Under the previous assumptions, the normalization map
 $S_{\overline{\eta}}'\to S_{\overline{\eta}}$ is a 
 homeomorphism, i.e., $S_{\overline{\eta}}$ has only 
 unibranch singularities (``cusps'').
 Moreover, if $p\geq3$, then the arithmetic genus of every cusp of 
 $S_{\overline{\eta}}$ is divisible by $(p-1)/2$.
\end{Theorem}

If the generic fiber $S_\eta$ has arithmetic genus one, 
and the fiber is not smooth over $k(B)$, then the normalization
of $S_{\overline{\eta}}$ is $\PP^1$. 
Also, $S_{\overline{\eta}}$ can have only one 
singularity, which, by Tate's theorem, must be
a cusp of arithmetic genus one.
Since $(p-1)/2$ divides this genus if $p\geq3$,
we find $p=3$ as only solution. 
Thus:

\begin{Corollary}
 Let $f:S\to B$ be a fibration from a smooth surface 
 whose generic fiber $S_\eta$ is a curve 
 of arithmetic genus one, i.e., $h^1(S_\eta,\OO_{S_\eta})=1$. 
 Then 
 \begin{enumerate}
  \item either $S_\eta$ is smooth over $k(B)$,
  \item or $S_{\overline{\eta}}$ is a singular rational curve 
   with one cusp.
 \end{enumerate}
 The second case can happen in characteristic $2$ and $3$ only.
\end{Corollary}

\begin{Definition}
 If the generic fiber of a fibration $S\to B$ is a smooth curve of
 genus one, the fibration is called {\em elliptic}.
 If the generic fiber is a curve that is not smooth over $k(B)$,
 the fibration is called {\em quasi-elliptic}, which can exist
 in characteristic $2$ and $3$ only.
\end{Definition}

Some authors require elliptic fibrations
to have a section, which we do not.
The literature is not consistent.

We refer to \cite{bm2} and \cite[Exercises 7.5 and 7.6]{Badescu} 
for examples of quasi-elliptic fibrations and to
\cite{bm3} for more on the geometry of quasi-elliptic 
fibrations.
For results on quasi-elliptic fibrations
in characteristic $3$, see \cite{Lang Quasi}.

We note that quasi-elliptic surfaces are always
uniruled, but may not be birationally ruled and
refer to Section \ref{subsec:quasi-elliptic uniruled}, where 
we discuss this in greater detail.

We also note that the situation gets more complicated in higher dimensions:
Mori and Saito \cite{Mori} constructed Fano $3$-folds $X$ in 
characteristic $2$ together with
fibrations $X\to S$, whose
generic fibers are conics in $\PP^2_{k(S)}$ 
that become non-reduced over $\overline{k(S)}$.
Such fibrations are called {\em wild conic bundles}.

\subsection{Canonical bundle formula}
Let $S$ be a smooth surface and
$f:S\to B$ be an elliptic or quasi-elliptic fibration.
Since $B$ is smooth, we obtain a decomposition
$$
 R^1 f_\ast\OO_S \,=\, {\cal L}\,\oplus\,{\cal T}\,,
$$
where ${\cal L}$ is an invertible sheaf and ${\cal T}$
is a torsion sheaf on $B$.
In characteristic zero, the torsion sheaf
$\cal T$ is always trivial.

\begin{Definition}
 Let $b\in B$ a point of the support of $\cal T$.
 Then, the fiber of $f$ above $b$ is called a
 {\em wild fiber} or, an {\em exceptional fiber}.
\end{Definition}

\begin{Proposition}
 Let $f:S\to B$ be a (quasi-)elliptic fibration,
 $b\in B$ and $F_b$ the fiber above $b$.
 Then the following are equivalent:
 \begin{enumerate}
  \item $b\in{\rm Supp}({\cal T})$, i.e., $F_b$ is a wild fiber,
  \item $h^1(F_b,\OO_{F_b})\geq2$,
  \item $h^0(F_b,\OO_{F_b})\geq2$.
 \end{enumerate}
 In particular, wild fibers are multiple fibers.
 Moreover, if $F_b$ is a wild fiber, then
 its multiplicity is divisible by $p$, and
 we have $h^1(S,\OO_S)\geq1$.
\end{Proposition}

The canonical bundle formula for relatively minimal (quasi-)elliptic fibrations
has been proved in \cite{bm2} - as usual, {\em relatively minimal}
means that there are no exceptional $(-1)$-curves in the fibers of
the fibration.

\begin{Theorem}[Canonical bundle formula]
 Let $f:S\to B$ be a relatively minimal 
 (quasi-)elliptic fibration from a smooth surface.
 Then,
 $$
  \omega_S\,\iso\,f^\ast(\omega_B\otimes{\cal L}^\vee)\otimes
  \OO_S\left(\sum_i a_i P_i\right),
 $$
 where 
 \begin{enumerate}
  \item $m_iP_i=F_i$ are the multiple fibers of $f$,
  \item $0\leq a_i< m_i$,
  \item $a_i=m_i-1$ if $F_i$ is not a wild fiber, and
  \item $\deg(\omega_S\otimes{\cal L}^\vee)=2g(B)-2+\chi(\OO_S)+{\rm length}({\cal T})$.
 \end{enumerate}
\end{Theorem}

For more results on the $a_i$'s we refer to
\cite[Proposition V.5.1.5]{Cossec; Dolgachev},
as well as to \cite{ku} for more details on wild fibers.

\subsection{Degenerate fibers of (quasi-)elliptic fibrations}
Usually, an elliptic fibration has fibers that are not smooth and
the possible cases have been classified by
Kodaira and N\'eron.
The list in positive characteristic is the same as in characteristic zero, cf.
\cite[Chapter V, \S1]{Cossec; Dolgachev} and
\cite[Theorem IV.8.2]{aec2}.
This is not such a surprise, as the
classification of degenerate fibers rests on
the adjunction formula and on matrices of
intersection numbers, and these numerics
do not depend on the characteristic of the ground
field.

Let us recall that the possible singular fibers together with
their Kodaira symbols are as follows (after reduction):
\begin{enumerate}
 \item An irreducible rational curve with a node as singularity (${\rm I}_1$).
 \item A cycle of $n\geq2$ rational curves (${\rm I}_n$).
 \item An irreducible rational curve with a cusp as singularity (${\rm II}$).
 \item A configuration of rational curves forming a root system of type
   $A_2^*$ (${\rm III}$), $\overline{A}_3$ (${\rm IV}$), 
   $\widetilde{E}_6$ (${\rm IV}^*$), $\widetilde{E}_7$ (${\rm III}^*$), 
   $\widetilde{E}_8$ (${\rm II}^*$) or $\widetilde{D}_n$ (${\rm I}^*_{n-4}$).
\end{enumerate}
In the first two cases, the reduction is called 
{\em multiplicative} or {\em semi-stable}, whereas in the last two
cases, it is called {\em additive} or {\em unstable}.
The latter names come from 
the theory of N\'eron models, see
\cite[Chapter IV]{aec2} or \cite{BLR}.
The former names are explained by the fact that
semi-stable reduction remains semi-stable after 
pull-back, whereas unstable reduction may become semi-stable
after pull-back.
In fact, for every fiber with unstable reduction there exists
a pull-back, whose reduction is semi-stable 
\cite[Proposition IV.10.3]{aec2}.

For an elliptic fibration $S\to B$ from a smooth surface, the
second Chern class (Euler number), can be expressed in terms of
the singular fibers by {\em Ogg's formula}
$$
  c_2(S) \,=\, \sum_i \nu(\Delta_i),
$$
where $i$ runs through the singular fibers, 
$\Delta_i$ is the {\em minimal discriminant} of the singular
fiber and $\nu$ denotes its valuation.
If a fiber has $n$ irreducible components, then this minimal discriminant
is as follows
$$
\nu(\Delta) \,=\,\left\{
\begin{array}{ll}
 1 \,+\, (n-1) & \mbox{if the reduction is multiplicative, i.e., of type } {\rm I}_n, \\ 
 2 \,+\, (n-1) \,+\,\delta & \mbox{if the reduction is additive}. 
\end{array}
\right.
$$
Here, $\delta$ is the {\em Swan conductor}, or, {\em wild part of the conductor}
of the fiber, which is zero if $p\neq2,3$.
We refer to \cite[Chapter IV, \S10]{aec2} for details and to
\cite[Proposition 5.1.6]{Cossec; Dolgachev} for a version for 
quasi-elliptic fibrations.

In a quasi-elliptic fibration, all fibers are additive and the
geometric generic fiber is of type ${\rm II}$, i.e., an 
irreducible rational curve with one cusp.
After an inseparable base change of the base curve $B'\to B$, the normalized
pull-back yields a fibration whose generic fiber is of genus zero,
see also Theorem \ref{quasi-elliptic base-change}.
All fibers are reduced or have multiplicity
equal to the characterstic $p=2,3$.
The list of possible geometric fibers is as follows
\cite[Corollary 5.2.4]{Cossec; Dolgachev}:
\begin{center}
\begin{tabular}{ccl}
   $p=3$  &:&  ${\rm II}$, ${\rm IV}$, ${\rm IV}^*$ and ${\rm II}^*$,\\
   $p=2$  &:&  ${\rm II}$, ${\rm III}$, ${\rm III}^*$, ${\rm II}^*$ 
   and ${\rm I}_n^*$.
\end{tabular}
\end{center}

Finally, we mention that if a (quasi-)elliptic fibration from a surface
has a section, then there exists a {\em Weierstra\ss\ model} 
\cite[Chapter 5, \S5]{Cossec; Dolgachev}, which is
more involved in characteristic $2,3$ than
in the other characteristics.

\section{Enriques--Kodaira classification}
\label{sec:Classification}

We now come to the Kodaira--Enriques classification of surfaces.
In positive characteristic, it is due to Bombieri and Mumford, see
\cite{Mumford}, \cite{bm2} and \cite{bm3}.
Let $S$ be a smooth projective surface of Kodaira dimension $\kappa(S)$.

\subsection{Negative Kodaira dimension}
First, let us recall and repeat Theorem \ref{birationally ruled}:

\begin{Theorem}
 If $\kappa(S)=-\infty$, then $S$ is birationally ruled.
\end{Theorem}

In fact, 
$\kappa(S)=-\infty$ is equivalent to $p_{12}(S)=0$, 
where $p_{12}$ is the $12$.th plurigenus 
\cite[Theorem 9.8]{Badescu}.
Moreover, although their minimal models are not unique, 
they have the same structure as in characteristic zero
by Theorem  \ref{thm:nagata}.
In Section \ref{uni and super}, we shall see
that uniruled surfaces in positive characteristic
may {\em not} fulfill $\kappa=-\infty$.

\subsection{Positive Kodaira dimension}
We recall that the {\em canonical ring} of a 
smooth and proper variety $X$ is defined to be
  $$
   R_{\rm can}(X) \,:=\, \bigoplus_{n\geq0} H^0\left(X,\,\omega_X^{\otimes n}\right)\,.
  $$
This said, we have the following fundamental result

\begin{Theorem}[Zariski--Mumford]
 The canonical ring $R_{\rm can}(S)$ of a smooth projective surface 
 is a finitely generated $k$-algebra.
 If $\kappa(S)\geq0$, then $R_{\rm can}(S)$ has transcendence
 degree $1+\kappa(S)$ over $k$.
\end{Theorem}

We refer to \cite{bm2} and \cite[Corollary 9.10]{Badescu}.
More generally, we refer to \cite[Chapter 14]{Badescu}
for a discussion of Zariski decompositions and finite generation
of the more general rings $R(S,D)$ for a $\QQ$-divisor $D$ on 
$S$.

For a surface with $\kappa(S)\geq1$ one studies the 
{\em Iitaka-fibration}
$$
  S \,\dashrightarrow\,{\rm Proj}\, R_{\rm can}(S)\,.
$$
By the theorem of Zariski--Mumford just mentioned, 
the right hand side is a projective variety of dimension $\kappa(S)$.

\begin{Theorem}
 Let $S$ be a minimal surface with $\kappa(S)=1$. 
 Then, (the Stein factorization of) the Iitaka fibration 
 is a morphism, 
 which is a relatively
 minimal elliptic or quasi-elliptic fibration.
\end{Theorem}

If $\kappa(S)=1$ and $p\neq2,3$, then the fibration is elliptic and
unique, $|mK_S|$ for $m\geq14$ defines this fibration, and $14$
is the optimal bound \cite{ku}.
The main difficulties dealt with by Katsura and Ueno \cite{ku}
are related to wild fibers.
We note that their bound $m\geq14$ is better than Iitaka's 
bound $m\geq86$ over the complex numbers, since over the complex numbers also
analytic surfaces that are not algebraic are taken into
account.

Although we will discuss surfaces with $\kappa(S)=2$, i.e., surfaces
of general type, in Section \ref{sec:GeneralType}, let us already anticipate
Theorem \ref{thm:canonical model}:

\begin{Theorem}
 Let $S$ be a minimal surface with $\kappa(S)=2$.
 Then, the Iitaka fibration 
 is a birational morphism that contracts all rational $(-2)$-curves
 and nothing more.
\end{Theorem}

\section{Kodaira dimension zero}
\label{sec:KodairaZero}

As in the complex case, surfaces in positive characteristic
that are of Kodaira dimension zero fall into four classes.
However, there are new subclasses of Enriques surfaces in characteristic $2$,
and new subclasses of hyperelliptic surfaces, so-called {\em quasi-hyperelliptic surfaces}, 
in characteristic $2$ and $3$.
In particular, there are no fundamentally new classes in characteristic $p\geq5$.

We start with a result that follows from
the classification, especially from the explicit classification of
(quasi-)hyperelliptic surfaces:

\begin{Theorem}
 Let $S$ be a minimal surface with $\kappa(S)=0$.
 Then, $\omega_S^{\otimes 12}\iso\OO_S$, and in particular,
 $p_{12}(S)=1$.
\end{Theorem}

The key to the classification of minimal surfaces with $\kappa(S)=0$ is to use $K_S^2=0$
to rewrite Noether's formula (see Section \ref{subsec:riemann-roch}) as follows:
$$
10\,+\,12p_g\,=\,8h^{0,1}\,+\,2\Delta\,+\,b_2,
$$
where $\Delta:=2h^{0,1}-b_1$ measures the defect of smoothness of 
$\Pic^0(S)$.
We have $0\leq\Delta\leq 2p_g$ in general, and $p_g\leq1$, since $\kappa(S)=0$.
Also, all terms in the above formula are non-negative, which gives
a finite list of possibilities, leading eventually to four classes,
see the introduction of \cite{bm2} or \cite[Chapter 10]{Badescu}.
Let us now discuss these classes in detail:

\subsection{Abelian surfaces}
These are two-dimensional Abelian varieties.
Their main invariants are as in characteristic zero:
$$
\begin{array}{ccccccc}
 \omega_S\iso\OO_S & \mbox{ } & p_g=1 & \mbox{ } & h^{0,1}=2 & \mbox{ } & h^{1,0}=2 \\
 \chi(\OO_S)=1 & & c_2=0 & & b_1=4 & & b_2=2 
\end{array}
$$
Abelian surfaces are usually studied within the framework
of Abelian varieties of arbitrary dimension.
There exists a huge amount of literature on Abelian 
varieties and their moduli spaces,
both in characteristic zero and in positive characteristic,
see, e.g., \cite{Mumford Abelian}.

By an (unpublished) result of Grothendieck
\cite[Theorem 5.23]{Illusie ExistenceTheorem}, 
Abelian varieties lift formally
to characteristic zero, see also Section \ref{subsec:general lift}.

For an Abelian variety $A$ of dimension $g$,
multiplication by $p$ is a finite
morphism. 
Its kernel $A[p]$ is a finite and flat group
scheme of length $p^{2g}$, and refer to
Section \ref{subsec:curves}, where we already discussed
elliptic curves ($g=1$).
The identity component $A[p]^0$ is infinitesimal
of length at least $g$. 
The  quotient $A[p]/A[p]^0$ is an \'etale group scheme
isomorphic to $(\ZZ/p\ZZ)^r$, for some $0\leq r\leq g$.
This quantity $r$ is called the {\em $p$-rank} of $A$.
For Abelian varieties of dimension at most two, the
$p$-rank can be detected by the Frobenius-action 
$F:H^1(A,\OO_A)\to H^1(A,\OO_A)$.

\begin{Definition}
 \label{def:Abelian surface ordinary} 
 An Abelian surface $A$ is called 
 \begin{enumerate}
   \item {\em ordinary} if $r=2$. Equivalently, $F$ acts bijectively
     on $H^1(A,\OO_A)$.
   \item {\em supersingular} if $r=0$. Equivalently, $F$ is zero on
     $H^1(A,\OO_A)$.
 \end{enumerate}
\end{Definition}

We remark that the image of the Albanese morphism of a uniruled surface
is at most one-dimensional.
Thus, an Abelian surface cannot be uniruled.
We note this in view of Shioda's notion of supersingularity and 
its connection with unirationality  discussed in 
Section \ref{subsec:Shioda supersingular}.

\subsection{K3 surfaces}
These surfaces have the following invariants:
$$
\begin{array}{ccccccc}
 \omega_S\iso\OO_S & \mbox{ } & p_g=1 & \mbox{ } & h^{0,1}=0 & \mbox{ } & h^{1,0}=0 \\
 \chi(\OO_S)=2 & & c_2=24 & & b_1=0 & & b_2=22 
\end{array}
$$
Their formal deformation spaces are smooth $W(k)$-algebras 
in any characteristic, i.e.,
the Bogomolov--Tian--Todorov
unobstructedness theorem for K3 surfaces in positive 
characteristic is true \cite{Rudakov Shafarevich}:

\begin{Theorem}[Rudakov--Shafarevich]
 \label{thm:rs}
 A K3 surface has no global vector fields.
 Thus,
 $$
  H^2(S,\Theta_S)\,\iso\, 
  H^2(S,\Omega^1_S\otimes\omega_S) \,\stackrel{{\rm SD}}{\iso}\,
  H^0(S,\Theta_S)^\vee\,=\,0
 $$
 where SD denotes Serre duality.
 In particular, deformations of K3 surfaces
 are unobstructed.
\end{Theorem}

For K3 surfaces over arbitrary fields, we have 
$h^2(\Theta_S)=h^{1,2}=h^{1,0}$ by Serre duality and 
$h^{0,1}=0$ by our list of invariants.
Over the complex numbers, vanishing of the former then 
follows easily from the Hodge symmetry $h^{1,0}=h^{0,1}$, 
which is induced by complex conjugation
and thus, may not hold over arbitrary ground fields.
The proof in positive characteristic of \cite{Rudakov Shafarevich}
makes heavy use of purely characteristic-$p$-techniques, see
Section  \ref{subsec:foliation surfaces}.
We note that over fields of positive characteristic and
in dimension three
the Bogomolov--Tian--Todorov unobstructedness 
theorem for Calabi--Yau varieties may fail, cf.
\cite{Hirokado CY} and \cite{Schroeer CY}.

The vanishing $H^2(S,\Theta_S)=0$ implies that K3 surfaces
possess formal lifts over the Witt ring.
Deligne \cite{Deligne K3} 
showed in fact (see Section \ref{sec:Witt} for more
on lifts):

\begin{Theorem}[Deligne]
 K3 surfaces lift projectively to characteristic zero.
\end{Theorem}

The moduli space of polarized K3 surfaces in positive and
mixed characteristic exists by a result of Rizov \cite{Rizov}.
However, it is still open, whether these moduli spaces
are irreducible.
What makes moduli spaces of (polarized) K3 surfaces
so difficult to come by, is that no local or global
Torelli theorems are known
(except for supersingular K3 surfaces, see Section \ref{subsec:K3 Artin}).

We come back to K3 surfaces in Section \ref{uni and super},
where we discuss arithmetic conjectures and conjectural 
characterizations of unirational K3 surfaces.

\subsection{Enriques surfaces}
In characteristic $p\neq2$ these surfaces have the following invariants:
$$
\begin{array}{ccccccc}
 \omega_S\not\iso\OO_S & \mbox{ } & \omega_S^{\otimes 2}\iso\OO_S & \mbox{ } & p_g=0 & \mbox{ } & h^{0,1}=0 \\
 \chi(\OO_S)=1 & & c_2=12 & & b_1=0 & & b_2=10 
\end{array}
$$
Moreover, the canonical sheaf $\omega_S$ defines an \'etale double cover
$\widetilde{S}\to S$, where $\widetilde{S}$ is a K3 surface.
Also, there always exist 
elliptic or quasi-elliptic fibrations.
Every such fibration has precisely two multiple fibers, both of
which are not wild. 

The most challenging case is characteristic $2$, where
Enriques surfaces are characterized by (here, $\equiv$ denotes
numerical equivalence) 
$$
\begin{array}{ccccccccc}
 \omega_S\equiv\OO_S & \mbox{ } &
 \chi(\OO_S)=1 & \mbox{ } & c_2=12 & \mbox{ } & b_1=0 & \mbox{ } & b_2=10 
\end{array}
$$
It turns out that $p_g=h^{0,1}\leq1$, see \cite{bm3}.
Since $b_1=0$, we conclude that the Picard scheme of
an Enriques surface with $h^{0,1}=1$ is not smooth.
In this case, Frobenius induces a map 
$F:H^1(S,\OO_S)\to H^1(S,\OO_S)$,
which is either zero or a bijection. 
We thus obtain three possibilities:

\begin{Definition}
 An Enriques surface (in characteristic $2$) is called
 \begin{enumerate}
  \item {\em classical} if $h^{0,1}=p_g=0$, and
  \item {\em non-classical} if $h^{0,1}=p_g=1$.
     Such a surface is called
     \begin{enumerate}
       \item {\em ordinary} if Frobenius acts bijectively on
          $H^1(S,\OO_S)$, and
       \item {\em supersingular} if Frobenius is zero on
          $H^1(S,\OO_S)$.
     \end{enumerate}
 \end{enumerate}
\end{Definition}

All three types exist \cite{bm3}.
We note that the terminology is inspired by Abelian surfaces,
see Definition  \ref{def:Abelian surface ordinary}.

In any characteristic, every Enriques surface possesses elliptic 
or quasi-elliptic fibrations. 
Such a fibration always has multiple fibers.
Moreover, if $S$ is classical, then every (quasi-)elliptic fibration
has precisely two multiple fibers, both of multiplicity two and
neither of them is wild.
If $S$ non-classical, then there is only one multiple fiber,
which is wild with multiplicity two.
Finally, if $S$ is non-classical and ordinary it does not possess
quasi-elliptic fibrations.
We refer to \cite[Chapter V.7]{Cossec; Dolgachev} for details.

As explained in \cite[\S3]{bm3}, every Enriques surface
possesses a finite and flat morphism of degree two
$$
\varphi\,:\,\widetilde{S}\,\to\,S
$$
such that $\omega_{\widetilde{S}}\iso\OO_{\widetilde{S}}$, 
$h^0(\widetilde{S},\OO_{\widetilde{S}})=h^2(\widetilde{S},\OO_{\widetilde{S}})=1$,
and $h^1(\widetilde{S},\OO_{\widetilde{S}})=0$, i.e., $\widetilde{S}$
is ``K3-like''.
More precisely, in characteristic $\neq2$, or, if $S$
is non-classical and ordinary, 
then $\varphi$ is \'etale of degree two and
$\widetilde{S}$ is in fact a K3 surface.
However, in the remaining cases in characteristic $2$, 
$\widetilde{S}$ is only an
integral Gorenstein surface that need not even be normal, since
$\varphi$ is a torsor under an infinitesimal group scheme.
In any case and any characteristic, $\widetilde{S}$ is birational to the complete
intersection of three quadrics in $\PP^5$, see
\cite{Liedtke Enriques}, which generalizes results of 
Cossec \cite{Cossec Enriques} and
Verra \cite{Verra}.

Moreover, the moduli space of Enriques surfaces in
characteristic $p\neq2$ is irreducible, unirational,
smooth and $10$-dimensional.
In characteristic $2$, it consists of two irreducible,
unirational, and $10$-dimensional components, whose intersection is
$9$-dimensional.
This intersection corresponds to non-classical supersingular
surfaces, and outside their intersection the $10$-dimensional
components parametrize non-classical ordinary, and 
classical Enriques surfaces, respectively.
We refer to \cite{Liedtke Enriques} for details, as well as
\cite{Ekedahl SB Enriques} for a complementary  approach.

We refer to \cite{bm3}, \cite{Lang 1}, \cite{Lang 2},
\cite{Liedtke Enriques}, and, 
of course, to \cite{Cossec; Dolgachev} for more details 
and partial classification results.

\subsection{(Quasi-)hyperelliptic surfaces}
\label{subsec:hyperelliptic surfaces}
In characteristic $p\neq2,3$, these surfaces have the following invariants:
$$
\begin{array}{ccccccc}
 \omega_S\not\iso\OO_S & \mbox{ } & \omega_S^{\otimes 12}\iso\OO_S & \mbox{ } & p_g=0 & \mbox{ } & h^{0,1}=1 \\
 \chi(\OO_S)=0 & & c_2=0 & & b_1=2 & & b_2=2 
\end{array}
$$
Moreover, these surfaces are equipped with two elliptic fibrations:
one is the Albanese fibration $S\to E$, where $E$ is an elliptic curve,
and then, there exists a second fibration $S\to\PP^1$.
It turns out that all these surfaces arise as quotients
$$
 S\,=\, (E\times F)/G\,,
$$
where $E$ and $F$ are elliptic curves, and $G$ is a group acting faithfully on
both, $E$ and $F$.
The quotient yielding $S$ is via the diagonal action.
In particular, the classical list 
of Bagnera--DeFranchis (see \cite[List 10.27]{Badescu}) gives all classes.

The more complicated classes arise in characteristic $2$ and $3$.
First, these surfaces have invariants
(again, $\equiv$ denotes numerical equivalence) 
$$
\begin{array}{ccccccccc}
 \omega_S\equiv\OO_S & \mbox{ } &
 \chi(\OO_S)=0 & \mbox{ } & c_2=0 & \mbox{ } & b_1=2 & \mbox{ } & b_2=2 
\end{array}
$$
It turns out that $1\leq p_g+1=h^{0,1}\leq2$, 
and that surfaces with $h^{0,1}=2$ are precisely those with non-smooth
Picard schemes \cite{bm3}.

In any case, the Albanese morphism $S\to\Alb(S)$ is onto
an elliptic curve, and its generic fiber is a curve of genus one.
This motivates the following

\begin{Definition}
 The surface is called {\em hyperelliptic}, if $S\to\Alb(S)$
 is an elliptic fibration, and {\em quasi-hyperelliptic} if this
 fibration is quasi-elliptic.
\end{Definition}

In both cases, there exists a second fibration $S\to\PP^1$,
which is always elliptic.
Finally, for every (quasi-)hyperelliptic surface $S$, there exists
\begin{enumerate}
 \item an elliptic curve $E$,
 \item a curve $C$ of arithmetic genus one, which is smooth if 
  $S$ is hyperelliptic, or rational with a cusp if $S$ is
  quasi-hyperelliptic, and
 \item a finite and flat group scheme $G$ 
   (possibly non-reduced),
   together with embeddings
   $G\to\Aut(C)$ and $G\to\Aut(E)$, where $G$ acts by translations
   on $E$,
\end{enumerate}
such that $S$ is isomorphic to 
$$
   S \,\iso\, (E\times C)/G\,.
$$
The Albanese map arises as projection onto $E/G$ with fiber $C$,
and the other fibration onto $C/G\iso\PP^1$ is elliptic with
fiber $E$.

It turns out that $G$ may contain infinitesimal subgroups,
which gives rise to new cases even for hyperelliptic surfaces.
In any case, it turns out that $G$ is Abelian.
This implies that the especially large and non-Abelian 
automorphism groups 
of elliptic curves with $j=0$ in characteristic $2$ and $3$ 
(see Section \ref{subsec:curves}) do {\em not} give rise to new classes. 
We refer to \cite{bm2} for the complete classification
of hyperelliptic surfaces and to \cite{bm3} for
the classification of quasi-hyperelliptic surfaces.

An interesting feature in characteristic $2$ and $3$ is the
possibility that $G$ acts trivially on $\omega_{E\times C}$,
and thus, the canonical sheaf on $S$ is trivial.
In this case, we find $p_g=1$, $h^{0,1}=2$, and the Picard
scheme of $S$ is not reduced.

\section {General type}
\label{sec:GeneralType}

In this section we discuss surfaces of general type, and 
refer to \cite[Chapter VII]{bhpv} and the references
given there for the corresponding results over the 
complex numbers.

\subsection{Pluricanonical maps}
\label{subsec:pluricanonical}
Let $S$ be a minimal surface of general type.
Clearly,
$$
    K_S^2\,>\,0
$$ 
since some pluricanonical map has a
two-dimensional image.
However, we shall see below that Castelnuovo's inequality
$c_2>0$ may fail.
Let us recall that a {\em rational $(-2)$-curve}
is a curve $C$ on a surface with $C\iso\PP^1_k$
and $C^2=-2$. 

\begin{Theorem}
 \label{thm:canonical model}
  Let $S$ be a minimal surface of general type.
 Then, the (a priori) rational Iitaka fibration to the canonical model
 $$
     S\,\dashrightarrow\,S_{\rm can}\,:= \Proj R_{\rm can}(S)\,=\,
     \Proj\,\bigoplus_{n\geq0} H^0(S,\,\omega_S^{\otimes n})
 $$
 is a birational morphism that contracts all rational $(-2)$-curves
 and nothing more.
\end{Theorem}

Rational $(-2)$-curves form configurations, whose intersection
matrices are negative definite.
These matrices are Cartan matrices and correspond to
Dynkin diagrams of type $A$, $D$ and $E$.
In particular, the morphism $S\to S_{\rm can}$ contracts these curves
to DuVal singularities (also known as canonical singularities, or
rational double points), 
see \cite[Chapter 3]{Badescu}, as well as Artin's original papers 
\cite{Artin Sing1} and \cite{Artin Sing2}.

Bombieri's results on pluricanonical systems were
extended to positive characteristic in \cite{ek},
and refined in \cite{sb} and \cite{cfhr}, and
we refer to these articles for more results.
We give a hint of how to modify the classical proofs
below.
Also, the reader who is puzzled by the possibility
of purely inseparably uniruled surfaces of general type
in the statements below might want to look at 
Section \ref{uni and super} first.

\begin{Theorem}[Ekedahl, Shepherd-Barron]
  Let $S$ be a minimal surface of general type and
  consider the linear system $|m K_S|$ on the 
  canonical model $S_{\rm can}$:
  \begin{enumerate}
   \item it is ample for $m\geq5$ or if $m=4$ and $K_S^2\geq2$ or $m=3$ and $K_S^2\geq3$,
   \item it is base-point free for $m\geq4$ or if $m=3$ and $K_S^2\geq2$,
   \item it is base-point free for $m=2$ if $K_S^2\geq5$ and $p\geq11$ or $p\geq3$ and
      $S$ is not uniruled, and
   \item it defines a birational morphism for $m=2$ if $K_S^2\geq10$, $S$ has no pencil
      of genus $2$ curves and $p\geq11$ or $p\geq5$ and $S$ is not uniruled.
  \end{enumerate}
\end{Theorem}

Next, we have the following version of Ramanujam-vanishing
(see \cite[Theorem II.1.6]{ek} for the complete statement):

\begin{Theorem}[Ekedahl]
  \label{Ekedahl vanishing}
  Let $S$ be a minimal surface of general type and let
  ${\cal L}$ be an invertible sheaf that is numerically 
  equivalent to $\omega_S^{\otimes i}$ for some $i\geq1$.
  Then,  $H^1(S,{\cal L}^{\vee})=0$ 
  except possibly for certain surfaces in characteristic $2$
  with $\chi(\OO_S)\leq1$.
\end{Theorem}

On the other hand, minimal surfaces of
general type with $H^1(S,\omega_S^{\vee})\neq0$ in
characteristic $2$ do exist 
\cite[Proposition I.2.14]{ek}. 

Bombieri's proof of the above results over the complex numbers 
is based on vanishing theorems $H^1(S,{\cal L})=0$ for certain 
more or less negative invertible sheaves.
However, these vanishing results may fail in positive
characteristic, see \cite{ray} or Section \ref{subsec:kodaira fails}.
Ekedahl \cite{ek} overcomes this problem as follows:
he considers an invertible sheaf ${\cal L}$ 
and its Frobenius-pullback 
$F^*({\cal L})\iso{\cal L}^{\otimes p}$
as group schemes over $S$.
Then, Frobenius induces a short exact sequence of
group schemes (for the flat topology on $S$)
\begin{equation}
\label{alphaL}
  0\,\to\,\alpha_{\cal L}\,\to\,{\cal L}\,
  \stackrel{F}{\to}\,F^*({\cal L})\,\to\,0\,.
\end{equation}
By definition, $\alpha_{\cal L}$ is the kernel
of $F$, see also the definition of 
the group scheme $\alpha_p$ in Section \ref{subsec:group schemes}.
This $\alpha_{\cal L}$ is an infinitesimal group scheme 
over $S$ and can be thought of as a possibly non-trivial
family of $\alpha_p$'s over $S$.

Now, if ${\cal L}^{\vee}$  is ample, then
$H^1(S,{\cal L}^{\otimes\nu})=0$ for $\nu\gg0$
(Serre vanishing, see \cite[Theorem III.7.6]{Hartshorne}).
In order to get vanishing of $H^1(S,{\cal L})$,
we assume that this is not the case and replace 
${\cal L}$ by some ${\cal L}^{\otimes \nu}$ such that 
$H^1(S,{\cal L})\neq 0$ and $H^1(S,{\cal L}^{\otimes p})=0$.
Then, the long exact sequence in cohomology for
(\ref{alphaL}) yields
$$
   H^1_{\rm fl}(S,\alpha_{\cal L})\neq0.
$$
Such a cohomology class corresponds to an 
$\alpha_{\cal L}$-torsor, which implies that there exists a
purely inseparable morphism of degree $p$
$$
    Y\,\stackrel{\pi}{\to}\,S,
$$
where $Y$ is an integral Gorenstein surface, whose
dualizing sheaf satisfies
$\omega_Y\iso\pi^*(\omega_S\otimes{\cal L}^{p-1})$.
(The subscript ${\rm fl}$ in the cohomology group above
denotes the flat topology, which is needed since $\alpha_{\cal L}$-torsors 
are usually only locally trivial with respect to the flat topology.)

For example, suppose $S$ is of general type and 
${\cal L}=\omega_S^{\otimes(-m)}$ for some 
$m\geq1$.
Then {\em either} $H^1(S,{\cal L})=0$ and one proceeds as
in the classical case
{\em or} there exists an inseparable cover
$Y\to S$, where 
$\omega_Y^\vee$ is big and nef.
The second alternative implies that $S$ is inseparably
dominated by a surface of special type, namely $Y$,
and a further analysis 
of the situation leads {\em either} to a contradiction 
(establishing the desired vanishing result) {\em or} an explicit 
counter-example to a vanishing result.
For example, Theorem \ref{Ekedahl vanishing} is proved
this way.

\subsection{Castelnuovo's inequality}
\label{subsec:castelnuovo}
Over the complex numbers, surfaces of general type satisfy
{\em Castelnuovo's inequality} $c_2>0$.
In \cite{ray}, Raynaud constructs minimal surfaces of general type 
with $c_2<0$ in every characteristic $p\geq5$, i.e., this inequality
fails.
On the other, we have the following structure result:

\begin{Theorem}[Shepherd-Barron \cite{sb2}]
 \label{thm:c2 negative}
 Let $S$ be a minimal surface of general type.
 \begin{enumerate}
  \item If $c_2(S)=0$, then $S$ is inseparably dominated
    by a surface of special type.
  \item If $c_2(S)<0$, then the Albanese map $S\to\Alb(S)$
    has one-dimensional image, whose generic fiber is a singular
    rational curve.
    In particular, $S$ is uniruled.
 \end{enumerate}
 In characteristic $p\geq11$, surfaces of general type
 satisfy $\chi(\OO_S)>0$.
\end{Theorem}

We refer to \cite{sb2} for more detailed statements.
There do exist surfaces of general type with $c_2<0$,
but in view of Noether's formula $12\chi=c_1^2+c_2$,
one might ask whether the stronger inequality 
$\chi>0$ still holds for surfaces of general type
in positive characteristic.
At least, we have the following analog of a theorem of
Castelnuovo and DeFranchis:

\begin{Proposition}
 Let $S$ be a surface with $\chi(\OO_S)<0$. 
 Then,
 \begin{enumerate}
  \item $S$ is birationally ruled over a curve of genus $1-\chi(\OO_S)$, or
  \item $S$ is quasi-elliptic of Kodaira dimension $\kappa=1$ and 
    $p\leq3$, or
  \item $S$ is a surface of general type and $p\leq7$.
 \end{enumerate}
\end{Proposition}

\prf
If $\kappa=-\infty$, then $S$ is birationally ruled over a curve of genus
$1-\chi(\OO_S)$ and we get the first case.
Also, by the explicit classification, there are no surfaces with $\kappa=0$
and $\chi(\OO_S)<0$.
For $\kappa=2$, this is \cite[Theorem 8]{sb2}. 

If $\kappa=1$, then $S$ admits a (quasi-)elliptic fibration $S\to B$, say 
with generic fiber $F$.
Also, we may assume that $S$ is minimal.
In case $F$ is smooth then \cite{dolgachev} yields 
$c_2(S)\geq e(F)\cdot e(B)=0$, where $e$ denotes the Euler number.
Since $c_1^2(S)=0$ for a relatively minimal (quasi-)elliptic fibration,
Noether's formula yields $\chi(\OO_S)=0$.
Thus, if $\chi(\OO_S)<0$, then the fibration must be quasi-elliptic and such
surfaces exist for $p\leq3$ only.
\qed\medskip

Quasi-elliptic surfaces with $\chi(\OO_S)<0$ in characteristic $p\leq3$
can be found in \cite{ray}, i.e., the first two cases of the previous
proposition do exist.
On the other hand, it is still unknown whether there do
exist surfaces of general type with $\chi(\OO_S)\leq0$.

\subsection{Noether's inequality}
\label{subsec:Noether}
Every minimal surface of general type fulfills 
$$
   K_S^2\,\geq\,2p_g(S)\,-\,4 \mbox{ \qquad {\em (Noether's inequality)}.}
$$
Moreover, if the canonical map is composed with a pencil, then
$$
K_S^2\,\geq\,3p_g(S)\,-\,6\mbox{ \qquad {\em (Beauville's inequality)}}
$$
holds true.
If the canonical map is birational onto its image, then
$$
K_S^2\,\geq\,3p_g(S)\,-\,7\mbox{ \qquad {\em (Castelnuovo's inequality)}}
$$
holds true, see \cite{Liedtke Horikawa} and \cite{Liedtke Horikawa2}.
In particular, this area of geography of surfaces of general type
behaves as over the complex numbers.

We recall that surfaces that are extremal with respect to Noether's inequality
are called {\em Horikawa surfaces}.
More precisely, an {\em even Horikawa surface} is a minimal surface
of general type with $K^2=2p_g-4$, whereas an {\em odd
Horikawa surface} satisfies $K^2=2p_g-3$.
These surfaces are classified in arbitrary characteristic 
in \cite{Liedtke Horikawa} and \cite{Liedtke Horikawa2}.
Basically,
the same structure results as over the complex numbers hold for them:
most of them arise as double covers
of rational surfaces via their canonical maps.
In characteristic $2$, the canonical map may
become purely inseparable, and then, the corresponding
Horikawa surfaces are unirational, see also 
Section \ref{uni and super}.
We refer to the aforementioned articles for precise classification
results, description of the moduli spaces, as well as the description of the
subsets in these moduli spaces corresponding to surfaces with inseparable
canonical maps.
Finally, unirational Horikawa surfaces in characteristic $p\geq3$
were systematically constructed in \cite{Liedtke; Schuett}.

Also, Beauville's result that minimal surfaces of general type with
$K^2<3p_g-7$ are double covers of rational surfaces via their
canonical maps still holds in positive characteristic
\cite{Liedtke Horikawa2}.

\subsection{Bogomolov--Miyaoka--Yau inequality}
\label{subsec:bmy}
A minimal surface of general type over the complex
numbers fulfills $K_S^2\leq9\chi(\OO_S)$ or, equivalently,
$K_S^2\leq3c_2(S)$.
This is proved using analytic methods from
differential geometry.
Moreover, by a theorem of Yau, surfaces with $c_1^2=3c_2$
are uniformized by the complex $2$-ball and thus,
these surfaces are rigid by a theorem of Siu. 

Minimal surfaces of general type with $c_2\leq0$ (counter-examples to
Castelnuovo's inequality) provide counter-examples to the
Bogomolov--Miyaoka--Yau inequality.
But even if $c_2$ is positive, it may fail, as shown by
Parshin \cite{Parshin} and Szpiro \cite[Section 3.4.1]{Szpiro}.
More precisely, they construct series of examples, where 
$c_2$ is bounded and $c_1^2$ tends to infinity.
Let us also mention the  counter-examples of 
\cite[Kapitel 3.5.J]{bhh},
where covers of $\PP^2$ ramified over special
line configurations that only exist in positive 
characteristic are used.
Similar constructions appeared in \cite{Easton}.

Since Parshin's counter-examples have highly non-reduced
Picard schemes, he asked in \cite{Parshin letter}, whether 
surfaces of general type with reduced Picard schemes 
satisfy the Bogomolov--Miyaoka--Yau inequality.
Also this turns out to be wrong by the examples of
Jang \cite{Jang}.

In \cite[Remark (i) to Proposition 2.14]{ek},
a $10$-dimensional family of surfaces with 
$K^2=9$ and $\chi(\OO_S)=1$ in characteristic $2$ 
is constructed, i.e., rigidity on the 
Bogomolov--Miyaoka--Yau line fails.

On the other hand, there is the following positive
result \cite{sb}

\begin{Theorem}[Shepherd-Barron]
 If $S$ is a minimal surface of general type in 
 characteristic $2$ that lifts over $W_2(k)$
 then $c_1^2(S)\leq4c_2(S)$ holds true.
\end{Theorem}

We refer to \cite{sb} for results circling around
Bogomolov's inequality $c_1^2({\cal E})\leq 4c_2({\cal E})$  
for stable rank $2$ vector bundles.

\subsection{Global vector fields}
The tangent space to the automorphism group scheme
of a smooth variety is isomorphic to the space
of global vector fields.
Since a surface of general type has only finitely many
automorphisms, this implies that there are no
global vector fields on a surface of general type
in characteristic zero.
However, in positive characteristic, the automorphism
group scheme of a surface of general type has still
finite length, but may contain infinitesimal subgroup schemes,
which have non-trivial tangent spaces.
Thus, infinitesimal automorphism group schemes
of surfaces of general type in positive characteristic 
give rise to non-trivial global vector fields.
For examples, we refer to Lang's article \cite{Lang vector fields}.

\subsection{Inseparable Albanese Morphisms}

If $X$ is smooth and projective, 
then the pullback of a nonzero global one-form from ${\rm Alb}(X)$ to $X$
via the Albanese morphism is again nonzero by a result of Igusa (see Section \ref{subsec:one-forms}).
Now, for a smooth projective curve, the Albanese map is trivial ($g=0$), an isomorphism ($g=1$), or 
a closed immersion ($g\geq2$).
Despite of Igusa's result, it is still possible that the Albanese morphism becomes purely
inseparable in positive characteristic:

\begin{Theorem}
 \label{thm:purely insep Albanese}
 Let $k$ be an algebraically closed field of positive characteristic $p$.
 Then, for every Abelian surface $A$ over $k$, there exist infinitely many surfaces of general type 
 (actually, with arbitrary large $K_S^2$, i.e., 
 these examples do not lie in some bounded families), whose Albanese morphisms are generically finite and purely inseparable
 onto $A$.
\end{Theorem}

{\sc Proof:}\,
Let $\mathcal L$ be a very ample line bundle on $A$ such that
$H^1(A,{\mathcal L}^{\otimes i})=0$ for all $-p<i<0$.
Then, for a generic section $s\in H^0(A,{\mathcal L}^{\otimes p})$, we consider the 
purely inseparable cover $\pi:S\to A$ of degree $p$ that is given by taking a $p$.th root of $s$.
That is, $\pi_\ast\OO_S\cong\bigoplus_{i=0}^{p-1}{\mathcal L}^{\otimes(-i)}$ and the $\OO_A$-algebra
structure is given by multiplication by $s:{\mathcal L}^{\otimes (-p)}\to\OO_A$.
Since $s$ is generic, $S$ has at worst canonical singularities of type $A_{p-1}$, see
\cite[Theorem 3.4]{Liedtke Horikawa2}. 
From $\omega_S\iso\pi^\ast(\omega_A\otimes{\mathcal L}^{\otimes(p-1)})$, we find that $\omega_S$ is ample, i.e., 
$S$ is the canonical model of a surface of general type.
Taking cohomology, and using the assumption $h^1({\mathcal L}^{\otimes i})=0$ for $-p<i<0$, we find
$h^1(\OO_S)=h^1(\OO_A)=2$.
Since $S$ has rational singularities, we find $h^1(\OO_{\widetilde{S}})=2$. 
This gives the estimate $b_1(\widetilde{S})\leq4$ for the first $\ell$-adic Betti number 
(we don't know yet whether $\Pic(\widetilde{S})$ is reduced).
On the other hand, we have a surjective morphism
$$
     \alpha\,:\,\widetilde{S}\,\to\, S\,\to\, A,
$$
onto an Abelian surface, giving the estimate $b_1(\widetilde{S})\geq4$. 
Thus, $b_1(\widetilde{S})=4$ (and it follows that $\Pic(\widetilde{S})$ is reduced, but we will not need this fact). 
In any case, the Albanese variety $A':={\rm Alb}(\widetilde{S})$ of $\widetilde{S}$ is an Abelian surface.
We claim that $\alpha$ is the Albanese morphism of $\widetilde{S}$:
 namely, let ${\rm alb}:\widetilde{S}\to A'$ be the Albanese map.
 By universality, $\alpha$ factors over ${\rm alb}$ showing that ${\rm alb}$ is generically finite onto $A'$.
 Also, the induced map $A'\to A$ is finite.
Since $\alpha$ is a generically finite of degree $p$, which is a prime, we get simply for degree
reasons, that either $\widetilde{S}$ is birational to an Abelian variety (which is not the case), or
$A$ and $A'$ are isomorphic.
Thus, $\alpha$ is the Albanese morphism.

Finally we note that every sufficiently ample line bundle on $A$ satisfies the assumptions of our construction.
In particular, for every $A$ we obtain infinitely many such $S$ with unbounded $K_S^2$.
\qed

\subsection{Non-classical Godeaux surfaces}
Since $K_S^2>0$ for a minimal surface of general type, it is
natural to classify surfaces with $K^2_S=1$.
It turns out that these fulfill $1\leq\chi(\OO_S)\leq3$
and thus, the lowest invariants possible are as follows:

\begin{Definition}
 A {\em numerical Godeaux surface} is a minimal surface of
 general type with $\chi(\OO_S)=K_S^2=1$.
 Such a surface is called {\em classical} if $p_g=h^{0,1}=0$
 and otherwise {\em non-classical}.
\end{Definition}

In characteristic zero or in characteristic $p\geq7$,
numerical Godeaux surfaces are classical
\cite{Liedtke Godeaux}.
Moreover, quotients of a quintic surface in $\PP^3$ 
by a free $\ZZ/5\ZZ$-action 
(this construction is due to Godeaux) provide
examples of classical Godeaux surfaces in characteristic 
$p\neq5$.
Classical and non-classical Godeaux surfaces in characteristic
$p=5$ have been constructed by Lang \cite{Lang Godeaux} and
Miranda \cite{Miranda Godeaux}.
Non-classical Godeaux surfaces in characteristic $p=5$
have been completely classified in \cite{Liedtke Godeaux} 
- it turns out that all of them arise as quotients of 
(possibly highly singular) quintic surfaces in $\PP^3$
by $\ZZ/5\ZZ$ or $\alpha_5$.
We finally note that non-classical Godeaux surfaces
are precisely those numerical Godeaux surfaces that have 
non-reduced Picard schemes.

Quite generally, 
for every $n$ there exists an integer $P(n)$
such that minimal surfaces of general type with
$K^2\leq n$ in characteristic $p\geq P(n)$ 
have a reduced Picard scheme \cite{Liedtke Picard}.
Thus, $P(1)=7$, but $P(n)$ is unbounded as a function of
$n$.

\subsection{Surfaces with $\mathbf{p_g=0}$}
For a minimal surface of general type with $p_g=0$ over the complex numbers,
the inequality $\chi(\OO_S)>0$ forces $h^{0,1}=0$,
thus $\chi(\OO_S)=1$, and then, the Bogomolov--Miyaoka--Yau
inequality implies $1\leq K^2\leq 9$.
Interestingly, these (in-)equalities hold over any field:

\begin{Proposition}
 Let $S$ be a minimal surface of general type with $p_g(S)=0$.
 Then, the equalities and inequalities
 $$
   h^{0,1}(S)=0,\mbox{ \quad }\chi(\OO_S)=1,\mbox{ \quad and \quad }1\leq K_S^2\leq9
 $$
 hold true.
\end{Proposition}

\prf
From $p_g=0$ we get $\chi(\OO_S)\leq1$.
Suppose first, that $\chi(\OO_S)=1$ holds.
Then $h^{0,1}=0$ and we find 
$b_1=0$, which yields $c_2=2-2b_1+b_2\geq3$.
But then, Noether's formula yields
$$
K_S^2\,=\,12\chi(\OO_S)-c_2(S)\,\leq\, 12-3\,=\,9,
$$
which gives the desired (in-)equalities.

If $\chi(\OO_S)\leq0$, then Noether's formula implies $c_2(S)<0$.
By Theorem \ref{thm:c2 negative},  the Albanese morphism
$S\to{\rm Alb}(S)$ is a fibration over a curve.
Thus, $b_2\geq\rho(S)\geq2$ using Igusa's inequality.
Next, we have $h^{0,1}=1-\chi(\OO_S)$ and in particular,
$b_1\leq2(1-\chi(\OO_S))$.
Thus,
$$
c_2(S)\,=\,2-2b_1+b_2\,\geq\, 2\,-\,4(1-\chi(\OO_S))+2 \,=\, 4\chi(\OO_S).
$$
But then, Noether's formula implies
$$
12\chi(\OO_S)\,=\,K_S^2\,+\,c_2(S)\,\geq\,4\chi(\OO_S)\,+\,K_S^2
$$
and we obtain
$$
0\,\geq\,8\chi(\OO_S)\,\geq\, K_S^2,
$$
a contradiction. 
Thus, the $\chi(\OO_S)\leq0$-case cannot happen and we are done.
\qed\medskip

The first examples of algebraically simply connected surfaces of general type
with $p_g=0$ were constructed by Lee and Nakayama \cite{lee nakayama}, 
by adapting $\QQ$-Gorenstein smoothing techniques to positive
characteristic.

\begin{Theorem}[Lee--Nakayama]
 There do exist algebraically simply connected surfaces of general type with
 $p_g=0$, all values $1\leq K^2\leq4$, and in all characteristics
 $p\geq3$.
\end{Theorem}

\addtocontents{toc}{{\medskip \bf Special Topics in Positive Characteristic}}

\section{Unirationality, supersingularity, finite fields, and arithmetic}
\label{uni and super}

In this and the following section we discuss more 
specialized characteristic-$p$ topics.
In this section, we circle around rationality,
unirationality, their effect on N\'eron--Severi
groups, and the formal Brauer group.
We discuss these for K3 surfaces, and
surfaces over finite fields.
Finally, we discuss zeta functions and the Tate conjecture.

\subsection{An instructive computation}
\label{subsec:instructive example}
To start with, let $\varphi:X\dashrightarrow Y$ be
a dominant and generically finite morphism in characteristic zero.
Then, the pull-back of a non-zero pluricanonical form is again
a non-zero pluricanonical form.
Thus, if $\kappa(X)=-\infty$, also $\kappa(Y)=-\infty$ 
holds true.
However, over a field of positive characteristic $p$, the example
$$
 \varphi:t\mapsto t^p,\mbox{ \quad and then \quad }
 \varphi^*(dt)\,=\,dt^p\,=\,pt^{p-1}dt\,=\,0
$$
shows that the pull-back of a non-zero pluricanonical
form may become zero after pullback.
In particular, the previous characteristic zero argument,
which shows that
the Kodaira dimension cannot increase under generically
finite morphisms, breaks down.

However, if $S$ is separably uniruled, i.e., if there
exists a dominant rational map $\PP^1\times C\dashrightarrow S$
such that the finite field extension $k(S)\subset k(\PP^1\times C)$ is
separable, then this phenomenon does not occur, 
we find $\kappa(S)=-\infty$ and applying 
Theorem \ref{birationally ruled}, we conclude

\begin{Theorem}
 A separably uniruled surface is birationally ruled.
\end{Theorem}

In particular, if a surface $S$ is separably unirational, then $p_2(S)=0$.
But being dominated by a rational surface, its Albanese map
is trivial and so $b_1(S)=0$.
Thus, $S$ is rational by Theorem \ref{thm:castelnuovo zariski},
and we have shown:

\begin{Theorem}
  A separably unirational surface is rational.
\end{Theorem}

\subsection{Zariski surfaces}
On the other hand, Zariski \cite{Zariski} gave the first examples
of (inseparably) unirational surfaces in positive characteristic that are not
rational:
for a generic choice of a polynomial 
$f(x,y)\in k[x,y]$ of sufficiently large degree,
\begin{equation}
\label{zariski}
   z^p\,-\,f(x,y)\,=\,0
\end{equation}
extends to an inseparable cover $X\to\PP^2$, where
$X$ has ``mild'' singularities and where usually $\kappa(\widetilde{X})\geq0$
for some resolution of singularities $\widetilde{X}\to X$.
By construction, we have an inclusion of function fields
$$
k(x,y)\,\subset\,k(\widetilde{X})\,=\,k(x,y)[\sqrt[p]{f(x,y)}]\,\subset\,k(\sqrt[p]{x},\sqrt[p]{y})
$$
i.e., $\widetilde{X}$ is unirational.
Surfaces that arise as desingularizations of covers of the
form (\ref{zariski}) are called {\em Zariski surfaces}.

\begin{Theorem}[Zariski]
 \label{thm:zariski surfaces}
 In every positive characteristic there do exist unirational surfaces
 that are not rational.
\end{Theorem}

However, we have seen in Theorem \ref{thm:castelnuovo zariski}
that rational surfaces are still characterized as those surfaces
that satisfy $h^{0,1}=p_2=0$.

\subsection{Quasi-elliptic surfaces}
\label{subsec:quasi-elliptic uniruled} 
If $S\to B$ is a quasi-elliptic 
fibration from a surface $S$
with generic fiber $F$, then there exists 
a purely inseparable extension $L/k(B)$ of degree 
$p={\rm char}(k)$, such that 
$F_L:=F\times_{\Spec k(B)} \Spec L$ is not normal, i.e., the cusp
``appears'' over $L$, see \cite{bm3}.
Thus, the normalization of $F_L$ is isomorphic to 
$\PP^1_L$, and 
we get the following result

\begin{Theorem} 
 \label{quasi-elliptic base-change}
 Let $S$ be a surface and $S\to B$ be a quasi-elliptic 
 fibration.
 Then, there exists a purely inseparable and dominant
 rational map $B\times\PP^1\dashrightarrow S$, i.e.,
 $S$ is (purely inseparably) uniruled.
\end{Theorem}

In particular, if $S\to\PP^1$ is a quasi-elliptic
fibration, then $S$ is a Zariski surface, and thus,
unirational.

\subsection{Fermat surfaces}
\label{subsec:fermat}
If the characteristic $p={\rm char}(k)$ does not divide 
$n$, then the {\em Fermat surface} $S_n$, i.e., the hypersurface
$$
   S_n \,:=\, \{\, x_0^n\,+\,x_1^n\,+\,x_2^n\,+\,x_3^n\,=\,0\,\}\,\subset\,\PP^3_k
$$
is smooth over $k$.
For $n\leq3$ it is rational, for $n=4$ it is K3, and for $n\geq5$ it is 
of general type.
Shioda and Katsura have shown in \cite{Shioda 1974}
and \cite{Shioda; Katsura} that

\begin{Theorem}[Katsura--Shioda]
 \label{thm:katsura shioda}
 For $n\geq4$ and $p\nmid n$, the Fermat surface $S_n$ 
 in characteristic $p$ is unirational  
 if and only if there exists a 
 $\nu\in\NN$ such that $p^\nu\equiv-1\mod n$. 
\end{Theorem}

Shioda  \cite{Shioda 1986} 
generalized this result to {\em Delsarte surfaces}.
The example of Fermat surfaces shows that
being unirational is quite subtle.
Namely, one can show that the generic hypersurface of 
degree $n\geq4$ in $\PP^3_k$ is {\em not} unirational, 
and thus, being unirational is not a deformation 
invariant.

From the point of view of Mori theory it is interesting to note that unirational
surfaces that are not rational are covered by {\em singular} rational curves.
However, (unlike in characteristic zero) it is not possible to smoothen 
these families -- after all, possessing a pencil of smooth rational curves
implies that the surface in question is rational.

\subsection{Fundamental group}
There do exist geometric obstructions to
unirationality:
being dominated by a rational surface, the Albanese morphism
of a unirational surface is trivial, and we conclude $b_1=0$.
Moreover, 
Serre \cite{Serre pi1} showed that the fundamental group
of a unirational surface is finite, and 
Crew \cite{Crew pi1} that it does {\em not} contain $p$-torsion
in characteristic $p$.
A subtle invariant is the formal Brauer group 
(see Section \ref{subsec:K3 Artin} below),
whose height can prevent a surface from being unirational (and that 
may actually be the only obstruction to unirationality for K3 surfaces).

\subsection{Horikawa surfaces}
\label{subsec:unirational Noether}
Let us recall from Section \ref{subsec:Noether}
that a minimal surface of general type is called
an {\em even Horikawa surface} if it satisfies 
$ K^2=2p_g-4$.
This unbounded class is particularly easy to handle, since all
such surfaces arise as double covers of rational surfaces.
In view of the previous paragraph, let us also mention that they
are algebraically simply connected.
In \cite{Liedtke; Schuett}, we constructed unirational
Horikawa surfaces in arbitrarily large characteristics
and for arbitrarily large $p_g$.
Thus, although the generic Horikawa surface is not unirational,
being unirational is nevertheless a common phenomenon.

\subsection{K3 surfaces and Shioda-supersingularity}
\label{subsec:Shioda supersingular}
We recall that the {\em Kummer surface} ${\rm Km}(A)$ of an Abelian
surface $A$ is the minimal desingularization of the quotient
of $A$ by the sign involution.
In characteristic $p\neq2$ the Kummer surface is always a K3 surface.
Shioda \cite{Shioda 1977} determined when such surfaces 
are unirational -
in particular, his result establishes the existence of unirational
K3 surfaces in every characteristic $p\geq3$:

\begin{Theorem}[Shioda]
 \label{thm:Kummer unirational}
 Let $A$ be an Abelian surface in characteristic $p\geq3$.
 Then, the Kummer surface ${\rm Km}(A)$ is unirational if and only
 if $A$ is a supersingular Abelian variety.
\end{Theorem}

We recall from Definition \ref{def:Abelian surface ordinary} 
that an Abelian variety is called supersingular
if its $p$-torsion subgroup scheme $A[p]$ is infinitesimal.

To explain the notion of supersingularity introduced
by Shioda \cite{Shioda 1974} let us recall from 
Section \ref{subsec:Igusa} that  Igusa's inequality states $\rho\leq b_2$,
where $\rho$ denotes the rank of the N\'eron--Severi group
and $b_2$ is the second Betti number.

\begin{Definition}
 \label{def:shioda-supersingular}
 A surface $S$ is called {\em supersingular in the sense of Shioda}
 if $\rho(S)=b_2(S)$ holds true.
\end{Definition}

This notion is motivated by the following result, also from
\cite{Shioda 1974}

\begin{Theorem}[Shioda]
  \label{thm:uniruled are supersingular}
  Uniruled surfaces are Shioda-supersingular.
\end{Theorem}

The unirationality results on Kummer and Fermat surfaces show
that these classes of surfaces are unirational if and only if 
they are supersingular in the sense of Shioda.
This leads to the following

\begin{Conjecture}[Shioda]
 \label{conj:k3 unirational}
 A K3 surface is unirational if and only if it is 
 Shioda-supersingular.
\end{Conjecture}

Apart from Kummer surfaces (basically Theorem \ref{thm:Kummer unirational}),
this conjecture is known to be true in characteristic $2$:
the N\'eron--Severi lattices of Shioda-supersingular
K3 surfaces have been classified in \cite{RS supersingular}
and using these results they show 

\begin{Theorem}[Rudakov--Shafarevich]
 Every Shioda-supersingular K3 surface in characteristic $2$
 possesses a quasi-elliptic fibration.
 In particular, these surfaces are Zariski surfaces and
 unirational.
\end{Theorem}

It is also is known in the following cases:
for supersingular K3 surfaces with Artin invariant $\sigma_0\leq6$ 
(see Definition   \ref{def:artin invariant})
and $p=3$ \cite{RS supersingular},
for $\sigma_0\leq3$ and $p=5$ \cite{Pho Shimada},
and for elliptic K3 surfaces with $p^n$-torsion section \cite{Ito Liedtke}.

\subsection{K3 surfaces and Artin-supersingularity}
\label{subsec:K3 Artin}
There exists yet another notion of supersingularity, apart from those of
Definition \ref{def:Abelian surface ordinary} and 
Definition \ref{def:shioda-supersingular}:
for a K3 surface $S$ over $k$, Artin \cite{Artin Supersingular}
considers
the functor that associates to every Artin-algebra $A$
over $k$ the Abelian group
$$
{\rm Br}\,:\,A\,\mapsto\,\ker\left(
H^2(S\times A, \OO_{S\times A}^\times)\,\to\,
H^2(S, \OO_{S}^\times)
\right)\,.
$$
This functor is pro-representable by a one-dimensional 
formal group law, the so-called {\em formal Brauer group} 
$\widehat{{\rm Br}}(S)$ of $S$.
(Of course, this functor can be studied for arbitrary
varieties, not just K3 surfaces.
Under suitable conditions, which are satisfied for K3 surfaces,
it it pro-representable by a formal group law of dimension
$h^{0,2}$, see \cite{Artin Mazur}.)

Over a field of characteristic zero, there exists
for every one-dimensional formal group
law an isomorphism (the logarithm)
to the additive formal group law $\widehat{\GG}_a$.
In positive characteristic, this need no longer be the
case, and every formal group law has a discrete invariant,
called the {\em height}.
The height $h$ is a strictly positive integer or infinity, and 
measures the complexity of multiplication by $p$ in the group law.
For example, $h=\infty$ means that 
multiplication by $p$ is equal to zero.
By a result of Lazard, $h$ determines 
the formal group law if the ground field is algebraically closed.
Over an algebraically closed field,
a one-dimensional group law of height $h=1$
is isomorphic to the multiplicative group law $\widehat{\GG}_m$,
whereas height $h=\infty$ corresponds to the additive
group law $\widehat{\GG}_a$.
We refer to \cite{Hazewinkel} for more on formal group laws
and to Artin's and Mazur's original article
\cite{Artin Mazur} for applications to geometry.

For K3 surfaces, the height $h$ of the formal
Brauer group satisfies $1\leq h\leq 10$ or $h=\infty$.
This follows from the fact that $b_2=22$ together with the formula
\begin{equation}
 \label{eqn:height}
  \rho(S)\,\leq\,b_2(S)\,-\,2h(S),
\end{equation}
which holds if $h\neq\infty$, see \cite{Artin Supersingular}.
Moreover, the height $h$ stratifies the moduli space 
of K3 surfaces:
K3 surfaces with $h=1$ -- these are called 
{\em ordinary} -- are open in families, and 
surfaces with $h\geq h_0+1$ form a closed subset inside
families of surfaces with $h\geq h_0$.
We refer to \cite{Artin Supersingular}, \cite{Ogus height}, and
\cite{Katsura Geer} for more on the geometry of the
height stratification of the moduli space.

\begin{Definition}
 A K3 surface is called {\em supersingular in the sense of
 Artin}, if its formal Brauer group has infinite height.
\end{Definition}

Shioda-supersingular K3 surfaces are Artin-supersingular,
which follows from Formula (\ref{eqn:height}).
To prove the converse direction, one first reduces
to the case of finite fields, where it follows from the
Tate conjecture for K3 surfaces with $h=\infty$ (see the discussion below).
For elliptic K3 surfaces with $h=\infty$, this latter conjecture was established by Artin
\cite{Artin Supersingular}, and for K3 surfaces possessing a degree $2$ polarization
by Rudakov, Shafarevich and Zink \cite{RSZ}. 
Finally, it was established by Charles \cite{Charles} and Maulik \cite{Maulik} for every
characteristic $p\geq5$.

\begin{Theorem}
 Let $X$ be a K3 surface in characteristic $p>0$.
 Assume that $p\geq5$, or that $X$ is elliptic, or that it possesses a
 degree $2$ polarization.
 Then, $X$ is Artin-supersingular if and only if it is Shioda-supersingular.
\end{Theorem}

To stratify the moduli space of Artin-supersingular K3 surfaces,
we consider their N\'eron--Severi groups.
The discriminant of the
intersection form on $\NS(S)$ of an Artin-supersingular K3 surface
$S$ is equal to
$$
{\rm disc}\,\NS(S)\,=\,\pm p^{2\sigma_0}
$$
for some integer $1\leq \sigma_0\leq 10$ by
\cite{Artin Supersingular}.

\begin{Definition}
  \label{def:artin invariant}
 The integer $\sigma_0$ is called the {\em Artin invariant}
 of the Artin-supersingular K3 surface.
\end{Definition}

In characteristic $p\geq3$,
Shioda-supersingular K3 surfaces with $\sigma_0\leq2$ are
Kummer surfaces of supersingular Abelian surfaces.
Their moduli space is one-dimensional but non-separated.
Moreover, there is precisely one such surface with
$\sigma_0=1$, and it arises as ${\rm Km}(E\times E)$,
where $E$ is a supersingular elliptic curve.
We refer to \cite{Shioda supersingular} and \cite{Ogus} for details,
and to \cite{Schroeer Kummer} for the description in
characteristic $2$.
From Theorem \ref{thm:Kummer unirational},
it follows that Shioda-supersingular K3 surfaces with
$\sigma_0\leq2$ are unirational (see also the
discussion at the end of Section \ref{subsec:Shioda supersingular}).

Ogus \cite{Ogus} established
a Torelli theorem for Shioda-supersingular K3 surfaces
with marked Picard lattices
in terms of crystalline cohomology.
Finally, we refer to \cite{RS K3} for further results on
K3 surfaces in positive characteristic.

\subsection{Zeta functions and Weil Conjectures}
If $X$ is a smooth and projective variety of dimension $d$ 
over a finite field $\FF_q$, then we can count 
the number $\#X(\FF_{q^n})$ of its $\FF_{q^n}$-rational points and form 
its {\em zeta function} :
$$
Z(X,t)\,:=\,\exp \left(\sum_{n=1}^\infty \#X(\FF_{q^n}) 
\frac{t^n}{n}\right)\,.
$$
Weil conjectured many properties of $Z(X,t)$,
and it was Grothendieck's insight that many of these properties
would follow from the existence of a suitable cohomology theory, namely
$\ell$-adic cohomology.
These conjectures are now known to hold by work of 
Deligne, Dwork, Grothendieck, Weil and others - 
we refer to \cite[Appendix C]{Hartshorne} for an overview and to
\cite{Milne} for details.
In particular, the zeta function is a rational
function of the form
$$
  Z(X,t)\,=\, \frac{P_1(S,t)\,\cdot\, P_3(S,t)\,\cdot\, ... \,\cdot\, P_{2d-1}(X,t)} {P_0(X,t)\,\cdot P_2(X,t)\,\cdot\,...\,\cdot\,P_{2d}(X,t)}\,,
$$
where each $P_i(X,t)$ is a polynomial with integral coefficients, with constant term $1$,
and of degree equal to the $i$.th Betti number $b_i(X)$.
In the extremal cases we have 
$P_0(X,t)=1-t$ and $P_{2d}(X,t)=1-q^dt$.
Moreover, over the complex numbers, these polynomials factor as
$$
  P_i(X,t)\,=\,\prod_{j=1}^{b_i(X)} (1\,-\,\alpha_{ij} t)\,,
$$
where the $\alpha_{ij}$ are complex numbers (in fact, algebraic integers) of absolute
value $q^{i/2}$.
Finally, there is a functional equation
$$
Z(X, \frac{1}{q^d\, t})\,=\,\pm\, q^{d E/2}\,t^E\,\cdot\,Z(X,t),
$$
where $E$ is the Euler number $c_d(X)=c_d(\Theta_X)$.

The Frobenius morphism $F_q:x\mapsto x^q$ acts trivially on $\FF_q$, and topologically generates
the absolute Galois group ${\rm Gal}(\overline{\FF}_q/\FF_q)$.
It also induces an $\FF_q$-linear morphism $F_q:X_{\overline{\FF}_q}\to X_{\overline{\FF}_q}$.
Now, an important characterization of $P_i(X,t)$ is that it
is equal to the characteristic polynomial $\det(1-t\cdot F_q^*)$ of the linear map $F_q^*$
induced by $F_q$ on $\Het{i}(X_{\overline{\FF}_q},\QQ_\ell)$. 
In fact, taking this as definition for the $P_i(X,t)$'s, and noting that $\FF_{q^n}$-rational
points of $X$ correspond to fixed-points under $F_q^{n}$, the rationality of the zeta-function
and the specific form of its factors as given above follow from Lefschetz fixed-point formulae for powers of 
$F_q^*$ on $\ell$-adic cohomology \cite[Appendix C.4]{Hartshorne}.

There is an injective Chern map $c_1:{\rm NS}(X_{\overline{\FF}_q})\otimes\QQ_\ell\to \Het{2}(X_{\overline{\FF}_q},\QQ_\ell)$,
which is equivariant with respect to the Galois-actions of ${\rm Gal}(\overline{\FF}_q/\FF_q)$ 
on both sides.
In particular, a non-torsion invertible sheaf ${\cal L}\in{\rm NS}(X)$ is Galois-invariant, and 
thus, $c_1({\cal L})$ is a non-trivial and Galois-invariant class in $\Het{2}(X_{\overline{\FF}_q},\QQ_\ell)$.
Since 
Frobenius topologically generates ${\rm Gal}(\overline{\FF}_q/\FF_q)$, 
we see that $c_1({\cal L})$ is an eigenvector of $F_q^*$ for the eigenvalue $q$.
In particular, $q^{-1}$ is a root of $P_2(X,t)$ and thus, $(1-qt)$ divides $P_2(X,t)$.
Applying this argument to the whole of ${\rm NS}(X)$, we find that 
$$
  (1-qt)^{\rho(X)}\mbox{ \quad divides \quad }P_2(X,t),
$$
where $\rho(X)$ denotes the rank of ${\rm NS}(X)$.
In \cite{Tate Conjecture}, Tate conjectured 
that the image $c_1({\rm NS}(X)\otimes\QQ_\ell)$ is not only a subspace, but
is in fact equal to the whole eigenspace of $F_q^*$ to the eigenvalue $q$.
We shall now discuss this conjecture in greater detail:

\subsection{Tate Conjecture}
\label{subsec:tate conjecture}
Let us now specialize to the case where $S:=X$ is a smooth and projective 
surface over $\FF_q$.
For the factorization of $Z(S,t)$, we have $P_0(S,t)=1-t$, and $P_4(S,t)=1-q^2t$.
Since $\Het{1}(X_{\overline{\FF}_q},\QQ_\ell)$ is Galois-equivariantly
isomorphic to $\Het{1}({\rm Alb}(X)_{\overline{\FF}_q},\QQ_\ell)$ via the Albanese morphism, 
we conclude $P_1(X,t)=P_1({\rm Alb}(X),t)$.
And finally, by Poincar\'e duality, we have $P_3(S,t)=P_1(S,qt)$.
Thus, the ``interesting'' part of $\ell$-adic cohomology and the zeta function is
encoded in $P_2(S,t)$.
We have also just seen that $(1-qt)^{\rho(S)}$ divides $P_2(S,t)$.

Now, suppose for a moment that $S_{\overline{\FF}_q}$ is Shioda-supersingular.
After possibly replacing $\FF_q$ by a finite extension, we may
assume that all divisor classes of $S_{\overline{\FF}_q}$ are defined
over $\FF_q$ and then, we have $P_2(S,t)=(1-qt)^{b_2(S)}$.
Moreover, if ${\rm Alb}(S)$ is trivial or a curve, then $Z(S,t)$  
is equal to the zeta function of a birationally ruled surface.
This fits perfectly to Conjecture \ref{conj:k3 unirational}.
Also, one might expect that if a surface over $\FF_q$
satisfies $P_2(S,t)=(1-qt)^{b_2(S)}$, then it is Shioda-supersingular.
This expectation would follow from the following, more general conjecture
of Tate \cite{Tate Conjecture}:

\begin{Conjecture}[Tate conjecture]
  \label{conj:tate}
  Let $S$ be a smooth and projective surface over $\FF_q$ and factor $P_2(S,t)$ as
  $$
    P_2(S,t)\,=\,\prod_{j=1}^{b_2(S)}(1-\alpha_{2,j}t),\mbox{ \quad with }\alpha_{2,j}\in\overline{\QQ}\cap\ZZ.
  $$
  Then, the N\'eron--Severi rank $\rho(S)$ is equal to the number of times $q$ occurs among the $\alpha_{2,j}$.
\end{Conjecture}

For an overview, we refer to \cite[Lecture 2]{Ulmer notes}.
For a relation of Tate's conjecture with Igusa's inequality and
a conjecture of Artin and Mazur on Frobenius eigenvalues 
on crystalline cohomology, we refer to
\cite[Remarque II.5.13]{Illusie deRham-Witt}.
Also, Artin and Tate \cite[(C)]{Artin Tate}
refined Tate's conjecture as follows:
let $D_1,...,D_\rho$ be independent classes
in $\NS(S)$ and set $B:=\sum_i \ZZ D_i$.
Let $\#{\rm Br}(S)$ be the order of the Brauer group, 
which is conjecturally finite.
Then

\begin{Conjecture}[Artin--Tate]
 \label{conj:artin-tate}
 We have
 $$
  P_2(S,q^{-s}) \,\lin\, 
  (-1)^{\rho(S)-1}\,\cdot\,
  \frac{\#{\rm Br}(S)\,\cdot\,
  \det( \{D_i\cdot D_j\}_{i,j})}
  {q^{\chi(\OO_S)-1+b_1(S)}\,\cdot\, (\NS(S):B)^2}
  \,\cdot\,(1-q^{1-s})^{\rho(S)}
 $$
 as $s$ tends to $1$
\end{Conjecture}

In fact, Conjectures \ref{conj:tate} and \ref{conj:artin-tate} are equivalent, as
shown up to $p$-power by Artin and Tate \cite{Artin Tate} and 
the full equivalence
was established by Milne \cite{Milne Tate}.

For elliptic surfaces, the Artin--Tate conjecture is a function field analog of
the Birch--Swinnerton-Dyer conjecture, see  \cite{Artin Tate} and \cite{Artin Swinnerton}.
For explicit examples, progress on this conjecture and interrelations,
we refer to \cite{Ulmer notes}.

The Tate conjecture is known in the following cases: 

\begin{Theorem}[Tate \cite{Tate Conjecture proof}]
  The conjectures of Tate and Artin--Tate hold for Abelian varieties and products of curves 
   over finite fields.
\end{Theorem}

Let us discuss what is known for K3 surfaces:
for elliptic K3 surfaces, it was established by Artin and Swinnerton-Dyer
\cite{Artin Swinnerton}.
For ordinary K3 surfaces, it was established by Nygaard \cite{Nygaard}, and for
K3's with finite height of the formal Brauer group 
and $p\geq5$ by Nygaard and Ogus \cite{Ogus; Nygaard}.
For K3 surfaces of infinite height (Artin-supersingular),
equipped with a polarization of degree $2$ it was established by Rudakov, Shafarevich and Zink \cite{RSZ},
if $p$ is large with respect to a polarization degree
by Maulik \cite{Maulik}, and for $p\geq5$ by Charles \cite{Charles}.
Thus, we obtain

\begin{Theorem}
  \label{tate conjecture}
  The conjectures of Tate and Artin--Tate hold for K3 surfaces in characteristic $p\geq5$.
\end{Theorem}

By  \cite{Lieblich}, this implies that
there exist only finitely many K3 surfaces defined over a fixed finite field of characteristic $p\geq5$.
This is similar to the situation for Abelian varieties:
by \cite{Zarhin}, there exist only finitely many Abelian varieties of a fixed
dimension over a fixed finite field.

Coming back to $P_2(S,t)$, we note that
Poincar\'e duality implies that if $\beta$ 
is among the $\alpha_{2,j}$, then so is $q/\beta$.
For K3 surfaces, using the fact that $\deg P_2(S,t)=b_2(S)=22$
is even, this has the following surprising consequence
(see \cite[Theorem 13]{bht} for a proof)

\begin{Theorem}[Swinnerton-Dyer]
 \label{K3 even rank}
 Let $S$ be a K3 surface over $\FF_q$, and
 assume that the Tate-conjecture holds for $S$.
 Then, the geometric N\'eron--Severi rank $\rho(S_{\overline{\FF}_q})$ 
 is even.
\end{Theorem}

Interestingly, there are more restrictions on $P_2(S,t)$
if $S$ is a K3 surface, than those coming from the Weil conjectures,
see \cite{Zarhin trans} and \cite{Elsenhans}.

Let us finally note that if we have an $\alpha_{2,j}$ in the factorization of $P_2(S,t)$
of some surface $S$ over $\FF_q$ that is not of the form $\mu\cdot q$, where
$\mu$ is a root of unity, then $S_{\overline{\FF}_q}$ is {\em not} Shioda-supersingular,
and thus, not unirational.
For example, the zeta function of a Fermat surface 
$S_n\subset\PP^3$ over $\FF_p$ can be computed
explicitly using Gau\ss- and Jacobi-sums.
From this, one concludes that if $(S_n)_{\overline{\FF}_p}$ is Shioda-supersingular, 
then there must exist a $\nu$ such that $p^\nu\equiv-1\mod n$,
see \cite{Shioda; Katsura} or Theorem \ref{thm:katsura shioda}.

\section{Inseparable morphisms and foliations}
\label{sec:Foliations}

In this section we study 
inseparable morphisms of height one in greater detail.
On the level of function fields this is Jacobson's
correspondence, a kind of Galois correspondence 
for purely inseparable field extensions.
However, this correspondence is not via automorphisms
but via derivations.
On the level of geometry, this translates into
$p$-closed foliations.
For surfaces, it simplifies to $p$-closed
vector fields.
For other overviews, we refer to
\cite{Ekedahl foliations} and \cite[Lecture III]{Miyaoka}.

\subsection{Jacobson's correspondence}
Let us recall the classical Galois correspondence:
given a field $K$ and a finite and {\em separable} extension $L$, there exists
a minimal Galois extension of $K$ containing $L$, the Galois closure
$K_{\rm gal}$ of $L$.
By definition, the Galois group $G={\rm Gal}(K_{\rm gal}/K)$ of this
extension is  the group of automorphism of $K_{\rm gal}$ over $K$, 
which is finite of degree equal to $[K_{\rm gal}:K]$.
Finally, there is a bijective correspondence between subgroups of $G$
and intermediate fields $K\subseteq M\subseteq K_{\rm gal}$. 
In particular, there are only finitely many fields between
$K$ and $K_{\rm gal}$.

In Section \ref{subsec:frobenius} we encountered
extensions of height one of a field $K$.
It turns out that automorphism of purely inseparable extensions
are trivial, and thus give no insight into these extensions.
However, there does exist a Galois-type correspondence for 
such extensions, {\em Jacobson's correspondence}
\cite[Chapter IV]{Jacobson III}.
Instead of automorphisms, one studies {\em derivations} over $K$:

Namely, let $L$ be a purely inseparable extension of height one of $K$,
i.e., $K\subseteq L\subseteq K^{p^{-1}}$, or, equivalently,
$L^p\subseteq K$.
We remark that $K^{p^{-1}}$ plays the role of a Galois closure of $L$.
Next, we consider the Abelian group
$$
  {\rm Der}(L) \,:=\,\{ \delta:K^{p^{-1}}\to K^{p^{-1}}\,,\,
  \mbox{ $\delta$ is a derivation and }\delta(L)=0 \}.
$$
Since $\delta(x^p)=p\cdot x^{p-1}\cdot \delta(x)=0$,
these derivations are automatically
$K$-linear and thus, ${\rm Der}(L)$ is a $K$-vector space.
Also, ${\rm Der}(L)$ is a subvector space of ${\rm Der}(K)$.
In case $K$ is of finite transcendence degree $n$ over some 
perfect field $k$, then ${\rm Der}(K)$ is $n$-dimensional.

Now, these vector spaces carry more structure:
if $\delta$ and $\eta$ are derivations, then in general their composition
$\delta\circ\eta$ is no derivation, which is why one studies their
{\em Lie bracket}, i.e., the commutator  $[\delta,\eta]=\delta\circ\eta-\eta\circ\delta$, 
which is again a derivation.
Now, over fields of positive characteristic $p$ it turns out that
the $p$-fold composite $\delta\circ...\circ\delta$ is again a derivation.
The reason is that expanding this composition the binomial coefficients
occurring that usually prevent this composition from being a derivation
are all divisible by $p$, i.e., vanish.
This $p$-power operation is denoted by $\delta\mapsto\delta^{[p]}$.
It turns out that the $K$-vector spaces ${\rm Der}(K)$ and
${\rm Der}(L)$ are closed under the
Lie bracket, as well as the $p$-power operation.

\begin{Definition}
  A {\em $p$-Lie algebra} or {\em restricted Lie algebra} is a Lie
  algebra over a field of characteristic $p$ together with a
  $p$-power map $\delta\mapsto\delta^{[p]}$ satisfying the
  axioms in \cite[Definition 4 of Chapter V.7]{Jacobson}.
\end{Definition}

We refer to \cite[Chapter V.7]{Jacobson} for general results on $p$-Lie algebras.

So far, we have associated to every finite and purely inseparable
extension $L/K$ of height one a sub-$p$-Lie algebra of ${\rm Der}(K)$.
Conversely, given such a Lie algebra $(V,{-}^{[p]})$, we may form
the fixed set
$$
   (K^{p^{-1}})^{  (V,{-}^{[p]}) }  \,:=\,
   \{ x\in K^{p^{-1}}\,|\,\delta(x)=0\,\,\forall\delta\in V\}\,,
$$
which is easily seen to be a field.
Since elements of $V$ are $K$-linear derivations, this field
contains $K$.
Moreover, by construction, it is contained in
$K^{p^{-1}}$, i.e., of height one.

\begin{Theorem}[Jacobson]
 There is a bijective correspondence
 $$
 \{ \mbox{ height one extensions of $K$ } \} 
   \leftrightarrow
 \{ \mbox{ sub-$p$-Lie algebras of ${\rm Der}(K)$ } \}\,.
 $$
\end{Theorem}

Let us mention one important difference to Galois theory:
suppose $K$ is of transcendence degree $n$ over an algebraically closed  
field $k$, e.g., the function field of an $n$-dimensional variety over $k$.
Then, the extension $K^{p^{-1}}/K$ is finite of degree $p^n$.
For $n\geq2$ there are  
infinitely many sub-$p$-Lie algebras of
${\rm Der}(K)$ and thus, {\em infinitely} many fields between 
$K$ and $K^{p^{-1}}$.

\subsection{Curves}
\label{curves comp Frob}
Let $C$ be a smooth projective curve over a 
perfect field $k$ with function field $K=k(C)$.
Then, the purely inseparable field extension $K^p\subset K$
is of degree $p$ and corresponds to the
$k$-linear Frobenius morphism $F:C\to C^{(p)}$.

Since every purely inseparable extension $L/K$ 
of degree $p$ is of the form $L=K[\sqrt[p]{x}]$
for some $x\in L$, such extensions are of height
one, i.e., $K\subseteq L\subseteq K^{p^{-1}}$.
Simply for degree reasons, we see that
the $k$-linear Frobenius morphism is the
only purely inseparable morphism of degree $p$ between 
normal curves.
Since every finite purely inseparable field extension 
can be factored successively into extensions of
degree $p$, we conclude

\begin{Proposition}
 \label{prop:frobenius curves}
 Let $C$ and $D$ be normal curves over a perfect field
 $k$ and let $\varphi:C\to D$ be a purely inseparable morphism
 of degree $p^n$.
 Then, $\varphi$ is the $n$-fold composite of the
 $k$-linear Frobenius morphism.
\end{Proposition}

\subsection{Foliations}
From dimension two on there are many more purely inseparable
morphisms than just compositions of Frobenius.
In fact, if $X$ is an $n$-dimensional
variety with $n\geq2$ over an algebraically closed field $k$,
then the $k$-linear Frobenius morphism has 
degree $p^n$ and it factors over infinitely many height 
one morphisms.

To classify height one morphisms $\varphi:X\to Y$
from a fixed smooth variety $X$ over a perfect field $k$,
we geometrize Jacobson's correspondence as follows:

\begin{Definition}
 A {\em ($p$-closed) foliation} on a smooth variety $X$ is a saturated
 subsheaf ${\cal E}$ of the tangent sheaf $\Theta_X$ that is closed
 under the Lie bracket (${\cal E}$ is {\em involutive}) and the $p$-power operation.
\end{Definition}

Then, Jacobson's correspondence translates into

\begin{Theorem}
 There is a bijective correspondence 
 $$
 \left\{ \begin{array}{l}
       \mbox{ finite morphisms $\varphi:X\to Y$ }\\
       \mbox{  of height one with $Y$ normal }
     \end{array} \right\} 
   \leftrightarrow
 \{ \mbox{ foliations in $\Theta_X$ } \}
 $$
\end{Theorem}

The saturation assumption is needed because an involutive and $p$-closed subsheaf
and its saturation (which will also be involutive and $p$-closed) 
define the same extension of function fields, and thus, correspond to the same normal variety.
We refer to \cite{Ekedahl foliations} or \cite[Lecture III]{Miyaoka}
for details.

Let us also mention \cite[Lecture III.2]{Miyaoka}, where a connection between
$p$-closed foliations and non-stability of tangent bundles, and uniruledness
of varieties (not only in positive characteristic, but also in characteristic zero!) is
discussed.

\subsection{Surfaces}
\label{subsec:foliation surfaces}
In order to
describe finite morphisms of height one $\varphi:X\to Y$
from a smooth surface onto a normal surface,
we have to consider foliations inside $\Theta_X$.
The sheaf $\Theta_X$ and its zero subsheaf correspond to
the $k$-linear Frobenius morphism and the identity,
respectively.
Thus, height one-morphisms of degree $p$
correspond to foliations of rank one inside $\Theta_X$.

To simplify our exposition, let us only consider $\Aff^2_k$, 
i.e., $X=\Spec R$ with $R=k[x,y]$ and assume that $k$ is
perfect.
Then, $\Theta_X$ corresponds to the $R$-module generated by
$\partial/\partial x$ and $\partial/\partial y$.
Now, a finite morphism of height one $\varphi:X\to Y$
with $Y$ normal corresponds to a ring extension
$$
  R^p=k[x^p,y^p]\,\subseteq\, S \,\subseteq\, R=k[x,y]\,,
$$
where $S$ is normal.
By Jacobson's correspondence,
giving $S$ is equivalent to giving a foliation
inside $\Theta_X$, which will be of rank one
if $S\neq R,R^p$.
This amounts to giving a regular vector field
$$
\delta \,=\, 
f(x,y)\,\frac{\partial}{\partial x}\,+\,
g(x,y)\,\frac{\partial}{\partial y}
$$
for some $f,g\in R$.
Since the Lie bracket of a $1$-dimensional Lie algebra is zero,
every rank one subsheaf of $\Theta_X$ is involutive.
Thus, we only have to check closedness under the $p$-power
operation, which translates into
$$
\delta^{[p]} \,=\, h(x,y) \cdot \delta\,\mbox{ \quad for some \quad }
h(x,y)\in R,
$$
i.e., $\delta$ is a {\em $p$-closed vector field}.

We may assume that $f$ and $g$ are coprime.
Then, the zero set of the ideal $(f,g)$ is of codimension two and
is called the  {\em singular locus} of the vector field.
It is not  difficult to see that $S$ is smooth over $k$ 
outside the singular locus of $\delta$, cf.
\cite{Rudakov Shafarevich}.

Finally, a purely inseparable morphism $\varphi:X\to Y$
is everywhere ramified, i.e., $\Omega_{X/Y}$ has support on
the whole of $X$.
Nevertheless, the canonical divisor classes of $X$ and $Y$
are related by a kind of Riemann--Hurwitz formula
and the role of the ramification divisor is played by
a divisor class that can be read off from the foliation,
see \cite{Rudakov Shafarevich}.

As an application, let us give the main result of \cite{Rudakov Shafarevich}:
let $S$ be a K3 surface over $k$, and suppose that we had
$H^0(S,\Theta_S)\neq0$.
Then, there exists in fact a 
$0\neq\delta\in H^0(S,\Theta_S)$ that is $p$-closed.
As explained above, this $\delta$ gives rise to an 
inseparable morphism $S\to S/\delta$.
A careful analysis of the hypothetical quotient $S/\delta$ and its
geometry finally leads to a contradiction, and we conclude
$H^0(S,\Theta_S)=0$, which proves 
Theorem \ref{thm:rs}.

\subsection{Quotients by group schemes}
\label{subsec:quot by alpha}
Let $X$ be a smooth but not necessarily proper
variety of any dimension over a perfect field $k$.
We have seen that a global section 
$0\neq\delta\in H^0(X,\Theta_X)$
gives rise to an inseparable morphism of degree $p$ and
height one if and only if $\delta$ is $p$-closed, i.e.,
$\delta^{[p]}=c\cdot\delta$ for some $c\in H^0(X,\OO_X)$.
Now, if $X$ is proper over $k$, then $c\in H^0(X,\OO_X)=k$, and after rescaling $\delta$,
we may in fact assume $c=1$ or $c=0$.

\begin{Definition}
 A vector field $\delta$ is called {\em multiplicative}
 if $\delta^{[p]}=\delta$ and it is called {\em additive}
 if $\delta^{[p]}=0$.
\end{Definition}

Let $\delta$ be additive or multiplicative.
Applying a (truncated) exponential series to $\delta$,
one obtains on $X$ an action of some finite and flat
group scheme $G$, which is infinitesimal of
length $p$,  see \cite[Section 1]{Schroeer Kummer}.
Then, the inseparable morphism 
$\varphi:X\to Y$ corresponding to $\delta$
is the quotient morphism $X\to X/G$.
Moreover, the $1$-dimensional $p$-Lie algebra
generated by $\delta$ is the $p$-Lie algebra of $G$,
i.e., the Zariski tangent space of $G$ 
with $p$-power map coming from Frobenius.
We recall from Theorem \ref{thm:oort tate} that
the only infinitesimal group schemes of length $p$
are $\alpha_p$ and $\mu_p$.
Putting these observations together, we obtain

\begin{Proposition}
 Additive (resp., multiplicative) vector fields correspond to 
 purely inseparable morphisms of degree $p$ that are quotients 
 by $\alpha_p$- (resp., $\mu_p$-) actions.
\end{Proposition}

This also explains the terminology for these vector fields:
$\alpha_p$ (resp., $\mu_p$)
is a subgroup scheme of the additive group
$\GG_a$ (resp., multiplicative group $\GG_m$).

\subsection{Singularities}
Let us finally assume that $X$ is a smooth surface and
let $\delta$ be a multiplicative vector field.
By \cite{Rudakov Shafarevich},
such a vector field can be written near a 
singularity in local coordinates $x$, $y$ as 
$$
 \delta \,=\, x\frac{\partial}{\partial x} \,+\,a\cdot y\frac{\partial}{\partial y}
\mbox{ \quad for some\quad }a\in\FF_p^\times\,.
$$

Let $\varphi:X\to Y$ be the inseparable morphism corresponding
to $\delta$.
In \cite{Hirokado Sing} it is shown that $Y$
has toric singularities of type $\frac{1}{p}(1,a)$.
Thus, quotients by $\mu_p$ behave very much like
cyclic quotient singularities in characteristic zero.
On the other hand, quotients by $\alpha_p$ are much more
complicated - the singularities need not even be rational
and we refer to \cite{Liedtke Uniruled} for examples.

\addtocontents{toc}{{\medskip \bf From Positive Characteristic to Characteristic Zero}}
\section{Witt vectors and lifting}
\label{sec:Witt}

This section deals with lifting to characteristic zero.
There are various notions of lifting, and the nicest ones are
are projective lifts over the Witt ring.
For example, in the latter case
Kodaira vanishing and degeneracy of the Fr\"olicher
spectral sequence hold true.
Unfortunately, although such lifts exist for curves,
they do not exist in general in dimension at least two.

\subsection{Witt vectors}
\label{subsec:Witt def}
Let $k$ be a field of positive characteristic $p$.
Moreover, assume that $k$ is perfect, 
e.g., algebraically closed or a finite field.

Then, one can ask whether there exist rings of 
characteristic zero having $k$ as residue field.
It turns out that there exists
a particularly nice ring $W(k)$,
the so-called {\em Witt ring}, or
{\em ring of Witt vectors},
which has the following properties:
\begin{enumerate}
\item $W(k)$ is a discrete valuation ring of characteristic zero,
\item the unique maximal ideal $\idealm$ of $W(k)$ 
 is generated by $p$ and the residue field 
 $R/\idealm$ is isomorphic to $k$,
\item $W(k)$ is complete with respect to the 
 $\idealm$-adic topology,
\item the Frobenius map $x\mapsto x^p$ on $k$ lifts to 
 a ring homomorphism of $W(k)$,
\item there exists an additive map $V:W(k)\to W(k)$,
 called {\em Verschiebung} (German for ``shift''), 
 which is zero on the residue field $k$ and
 such that multiplication by $p$ on $W(k)$ factors
 as $p=F\circ V=V\circ F$, and finally
\item every complete discrete valuation ring with
  quotient field of characteristic zero and
  residue field $k$ contains $W(k)$ as subring.
\end{enumerate}
We remark that the last property characterizes $W(k)$
up to isomorphism.

To obtain $W(k)$, one constructs successively 
rings $W_n(k)$, which are local
Artin rings of length $n$ with residue field $k$.
One has $W_1(k)=k$ and surjective projection maps
$W_{n+1}(k)\to W_n(k)$.
By definition, $W(k)$ is the projective limit over the 
$W_n(k)$, cf. \cite[Chapitre II.6]{Serre CorpsLocaux}.
The main example to bear in mind is the following:

\begin{Example}
 For the finite field $\FF_p$
 we have $W_n(\FF_p)\iso\ZZ/p^n\ZZ$ and thus,
 $$
     W(\FF_p)\,=\,\varprojlim\, \ZZ/p^n\ZZ
 $$
 is isomorphic to $\ZZ_p$, the ring of $p$-adic
 integers.
The maximal ideal of $W(\FF_p)$
 is generated by $p$ and $W(\FF_p)$ is complete with respect 
 to the $p$-adic topology.
 In this special case, $F$ is the identity on $W(\FF_p)$ and
 $V$ is multiplication by $p$.
\end{Example}

Witt's construction $W(-)$ makes sense
for every commutative ring $R$. 
However, already $W(k)$ for a non-perfect field $k$ is 
not Noetherian, and its maximal ideal is not generated
by $p$.
This is why we will assume $k$ to be perfect for the rest
of this section.
We refer to \cite[Chapitre II.6]{Serre CorpsLocaux} and
\cite{Hazewinkel} for more on Witt vectors.

\subsection{Lifting over the Witt ring}
\label{subsec:witt lift}
Let $X$ be a scheme of finite type over some perfect field $k$ 
of positive characteristic $p$.
Then, there are different notions of what it means to
{\em lift $X$ to characteristic zero}.
To make it precise, let $R$ be a ring of characteristic zero
with maximal ideal $\idealm$ and residue field
$R/\idealm\iso k$.
For example, we could have $R=W(k)$ and $\idealm=(p)$.

\begin{Definition}
 A {\em lift (resp. formal lift) of $X$ over $R$} is a  
 scheme (resp. formal scheme) ${\cal X}$ of finite type and 
 flat over $\Spec R$ (resp.
 ${\rm Spf}\: R$) with special fiber $X$.
\end{Definition}

In case $R=W(k)$, i.e., if $X$ admits
a (formal) lift over the Witt ring, many
``characteristic $p$ pathologies'' cannot
happen.
We have already encountered the following
results in Section \ref{sec:Cohomology}:
\begin{enumerate}
 \item if $X$ is of dimension $d\leq p$ and
   lifts over $W_2(k)$ 
   then its Fr\"olicher
   spectral sequence from Hodge to deRham-cohomology
   degenerates at $E_1$ by a result of Deligne and Illusie,
   see \cite{Deligne Illusie}
   and \cite[Corollary 5.6]{Illusie Frobenius},
 \item if $X$ is of dimension $d\leq p$ and
   lifts over $W_2(k)$,
   then ample line bundles satisfy Kodaira vanishing,
   see \cite{Deligne Illusie} and 
   \cite[Theorem 5.8]{Illusie Frobenius}, and
 \item if $X$ lifts over $W(k)$, then crystalline cohomology
   coincides with deRham-cohomology of ${\cal X}/W(k)$.
\end{enumerate}
Actually, the last property is the starting point of
crystalline cohomology, see the discussion in
Section \ref{subsec:crystalline}.
 
\begin{Example}
 Smooth curves and birationally ruled surfaces  
 lift over the Witt ring by
 Grothendieck's existence theorem
 \cite[Theorem 5.19]{Illusie ExistenceTheorem}.
\end{Example}

\subsection{Lifting over more general rings}
\label{subsec:general lift}
Let $R$ be an integral ring with maximal
ideal $\idealm$, residue field $R/\idealm\iso k$, and
quotient field $K$ of characteristic zero.
Let $X$ be a smooth projective variety over $k$,
let ${\cal X}$ be a lift of $X$ over
$\Spec R$, and denote its generic fiber by
${\cal X}_{K}\to\Spec K$.

After choosing a DVR dominating $(R,\idealm)$ and after passing to 
the $\idealm$-adic completion, we
may assume that $(R,\idealm)$ is a local
and $\idealm$-adically complete DVR.
By the universal property of the Witt ring, $R$ contains
$W(k)$ and $\idealm$ lies above
$(p)\subset W(k)$.
Thus, it makes sense to talk about the {\em ramification
index}, usually denoted by $e$, of $R$ over $W(k)$.
This ramification index 
is an absolute invariant of $R$.

To give a flavor of the subtleties that occur when dealing with
lifting problems, let us mention the following examples

\begin{enumerate}
 \item Abelian varieties admit formal lifts over the Witt ring by 
   an unpublished result of Grothendieck
   \cite[Theorem 5.23]{Illusie ExistenceTheorem}.
   However, to obtain algebraic lifts, one would like to
   have an ample line bundle on a formal lift in order to
   apply Grothendiecks' existence theorem, see 
   \cite[Theorem 4.10]{Illusie ExistenceTheorem}.
   However, even if one succeeds in doing so, this is usually at
   the prize that this new formal lift (which then is algebraic) 
   may exist over a {\em ramified} extension of the Witt ring only.
   For Abelian varieties, this was established by 
   Mumford \cite{Mumford Lifting}, and Norman and Oort \cite{Norman Oort}.

 \item K3 surfaces have unobstructed deformations
   by Theorem \ref{thm:rs}, and thus,
   admit formal lifts over the Witt ring.
   Deligne \cite{Deligne K3} has shown that one can lift with every 
   K3 surface also an ample line bundle, which gives an algebraic
   lifting - again at the prize that this lift may exist over
   ramified extensions of the Witt ring only.

\item By results of Lang \cite{Lang 1}, 
  Illusie \cite{Illusie deRham-Witt}, 
  Ekedahl and Shepherd-Barron \cite{Ekedahl SB Enriques}, and 
  \cite{Liedtke Enriques},
  Enriques surfaces - even in characteristic $2$ - lift to 
  characteristic  zero.
  However, the Fr\"olicher spectral sequence of a supersingular
  Enriques surface in characteristic $2$ does not degenerate 
  at $E_1$ by \cite[Proposition II.7.3.8]{Illusie deRham-Witt}.
  Thus, these latter surfaces only lift over 
  {\em ramified} extensions of the Witt ring, but not over
  the Witt ring itself.
\item Lang \cite{Lang 3} gave examples of hyperelliptic surfaces
  that lift to a ramified extension of $W(k)$ of ramification index $e=2$,
  but whose Fr\"olicher spectral sequences do not degenerate
  at $E_1$.
  Thus, these surfaces do not lift over $W(k)$.
  Rather subtle examples of non-liftable smooth elliptic fibrations
  were given by Partsch \cite{Partsch}.
\end{enumerate}

However, even if $X$ lifts ``only'' over a ramified extension of
the Witt ring, this does imply something:
flatness of ${\cal X}$ over $\Spec R$ implies that $\chi(\OO)$ of 
special and generic fiber coincide, and 
smoothness of ${\cal X}$ over $\Spec R$ implies that the $\ell$-adic
Betti numbers of special and generic fiber coincide.
For surfaces, we have additional results from
\cite[Section 9]{ku}:

\newcommand{\GenFib}{{{\cal S}_K}}

\begin{Theorem}[Katsura--Ueno]
 \label{katsura ueno theorem}
 Let ${\cal S}$ be a lift of the smooth projective surface
 $S$ over $\Spec R$ with generic
 fiber ${\cal S}_K$.
 Then,
 $$\begin{array}{ccccccc}
   b_i(S)&=&b_i(\GenFib)&\mbox{ \qquad }& c_2(S)&=&c_2(\GenFib)\\
   \chi(\OO_S)&=&\chi(\OO_\GenFib) & & K^2_S &=& K^2_\GenFib\\
   \kappa(S)&=&\kappa(\GenFib)
  \end{array}
 $$
 Moreover, $S$ is minimal if and only $\GenFib$ is minimal. 
\end{Theorem}

If $S$ is of general type then $P_n(S)=P_n(\GenFib)$ for $n\geq3$
since these numbers depend only on $\chi$ and $K^2$ by Riemann--Roch
and \cite[Theorem II.1.7]{ek}.
However, in general, $p_g(S)$ may differ from $p_g(\GenFib)$, as 
the examples in \cite{Serre Mexico} and \cite{Suh} show.
More precisely, Hodge invariants are semi-continuous, i.e., in general
we have 
$$
   h^{i,j}(S)\,\geq\,h^{i,j}(\GenFib)
   \mbox{ \quad for all \quad }i,j\,\geq0\,.
$$
In case of equality for all $i,j$, the Fr\"olicher spectral sequence 
of $S$ degenerates at $E_1$.
Theorem \ref{katsura ueno theorem} implies that from dimension two on
there exist smooth projective varieties
that do not admit any sort of lifting, namely:

\begin{Examples}
  Let $S$ be  
 \begin{enumerate}
   \item a minimal surface of general type with $K_S^2>9\chi(\OO_S)$, i.e.,
     violating the Bogomolov--Miyaoka--Yau inequality (see Section \ref{subsec:bmy}), or
   \item a quasi-elliptic surface with $\kappa(S)=1$ and
     $\chi(\OO_S)<0$ (see Section \ref{subsec:castelnuovo}).
 \end{enumerate}
  Then, $S$ does not admit an algebraic lifting whatsoever, i.e., 
   not even over a ramified extension of the Witt ring.
   The first example of such a smooth and projective variety
   that does not admit an algebraic lifting is due 
   to Serre  \cite{Serre non-lifting}.
\end{Examples}

For this and related questions, see also 
\cite[Section 5F]{Illusie ExistenceTheorem}.
Moreover, we have the following highly non-explicit result:
namely, ``Murphy's law'' holds for moduli spaces of surfaces of general
type with very ample canonical sheaves \cite{Vakil}.
Thus, we can find any kind of obstructed lifting behavior already
on surfaces, for example:

\begin{Theorem}[Vakil]
 For every integer $n>0$ and every prime $p>0$, there exists a smooth
 and projective surface over $\FF_p$ that lifts over
 $W_n(\FF_p)$ but not over $W_{n+1}(\FF_p)$.
\end{Theorem}

\subsection{Birational nature}
One can also ask to what extent liftability
is a birational invariant.
If $X$ and $Y$ are smooth, proper, and birational varieties of dimension at most $2$,
then their lifting behavior is the same.
However, in dimension $\geq3$, or when allowing canonical singularities
in dimension $2$, this is no longer the case.
We refer to \cite{Liedtke; Satriano} for details, some positive results, and (counter-)examples.

\subsection{Canonical lifts}
For an {\em ordinary} Abelian variety or K3 surface, there even exists a distinguished formal
lift over the Witt ring, the  {\em canonical lift}, or {\em Serre--Tate lift}.
Quite generally, ordinary means that Newton- and Hodge-polygons on crystalline cohomology coincide,
and we note that this property is open in equi-characteristic families.
For a $g$-dimensional Abelian variety $A$ over a field $k$ of characteristic $p$, 
being ordinary is equivalent to $A[p](\overline{k})\iso(\ZZ/p\ZZ)^g$, which is the maximum possible
(see also Definition \ref{def:Abelian surface ordinary}).
For a K3 surface $S$, being ordinary is equivalent to $h(\widehat{{\rm Br}}(S))=1$,
see Section \ref{subsec:K3 Artin}.
We refer to \cite{Messing} for details on canonical lifts of ordinary Abelian varieties.
For ordinary Abelian varieties, this canonical lift is characterized by the property
that the Frobenius morphism lifts. 
For the general case, we refer to \cite{Deligne canonical} and \cite{Katz}.
Finally, K3 surfaces with $h(\widehat{{\rm Br}}(S))<\infty$ still possess
{\em quasi-canonical lifts}, which has been used to prove the Tate conjecture for them, see
\cite{Ogus; Nygaard} and Section \ref{subsec:tate conjecture}.

\section{Rational curves on K3 surfaces}
\label{sec:rationalcurves}

In the final section we give an application of characteristic-$p$ and lifting 
techniques to a characteristic zero conjecture.
Namely, we show how infinitely many rational curves on 
complex projective K3 surfaces of odd Picard rank can be established
by reduction modulo $p$, then finding the desired
rational curves over finite fields, and eventually lifting cycles of them to characteristic zero.

\subsection{Rational curves}
Let $C$ be a smooth projective curve of genus $g$
over an algebraically closed field $k$.
Then, the Riemann--Hurwitz formula implies that 
if there exists a non-constant map $\PP^1\to C$ then
$g=0$, i.e., $C\iso\PP^1$.
Similarly, one can ask about non-constant maps from $\PP^1$ to higher
dimensional varieties, i.e., whether they exist and if so, how many,
whether they move in families, etc.
First of all, let us introduce the following notion 

\begin{Definition}
  A {\em rational curve} on a variety $X$ is a reduced and irreducible
  curve $C\subset X$ whose normalization is isomorphic to $\PP^1$.
\end{Definition}

Let us study rational curves on surfaces in detail:
clearly, if $S$ is a non-minimal surface, then every exceptional
$(-1)$-curve is a rational curve.

Also, since surfaces with $\kappa(S)=-\infty$ are birationally ruled
by Theorem \ref{birationally ruled}, they contain moving 
families of rational curves.
On the other extreme, 
Serge Lang \cite{Lang Diophant} conjectured that complex
surfaces of general type contain only finitely many rational
curves.
However, we note that uniruled surfaces of general type in positive
characteristic (see Section \ref{subsec:unirational Noether}) 
contain infinitely many rational curves.

\subsection{K3 surfaces}
In between these extremes lie surfaces of Kodaira dimension zero.
If $S$ is an Abelian variety, then every map $\PP^1\to S$
factors over the Albanese variety of $\PP^1$, which is a point.
Thus, Abelian varieties contain no rational curves at all.
On the other hand, there is the well-known

\begin{Conjecture}[Bogomolov]
 \label{infinity conjecture}
 A projective K3 surface contains infinitely many rational curves.
\end{Conjecture}

In characteristic zero, rational curves cannot move inside their linear
systems, for otherwise the K3 surface in question would have to be uniruled, 
which is impossible. 
But even in positive characteristic, where uniruled K3 surfaces do exist,
they are rather special, namely supersingular 
by Theorem \ref{thm:uniruled are supersingular}.
The first important step towards Bogomolov's conjecture is to establish
the existence of {\em at least one} rational curve, and we refer to
\cite{Mori Mukai} for the following result:

\begin{Theorem}[Bogomolov--Mumford]
  \label{one curve}
  Let $S$ be a projective K3 surface over an algebraically closed field,
  and let ${\cal L}$ be a non-trivial and effective invertible sheaf.
  Then, there exists a divisor $\sum_i n_i C_i$ inside $|{\cal L}|$,
  where $n_i\geq1$ and the $C_i$ are rational curves on $S$.
\end{Theorem}

For polarized K3 surfaces $(S,H)$, say, of degree $H^2=2d$,
there exists a moduli space ${\cal M}_{2d}$, which is smooth and
irreducible over the complex numbers, see, for example
\cite[Chapter VIII]{bhpv}.
Using degenerations of K3 surfaces to unions of rational surfaces, 
Chen \cite{Chen} showed, among other things,

\begin{Theorem}[Chen]
  A very general complex projective K3 surface in 
  ${\cal M}_{2d}$ contains infinitely many rational curves.
\end{Theorem}

Here, very general is meant in the sense that there exists a countable union
of analytic divisors inside ${\cal M}_{2d}$, outside of which the statement
is true.
Although this result strongly supports Conjecture \ref{infinity conjecture}, it 
does not give even a single example of a K3 surface containing infinitely many
rational curves!

\subsection{Explicit results}
It is shown in \cite[Section 4]{bt2}, or
\cite[Example 5]{bht} that complex projective
Kummer K3 surfaces contain infinitely many
rational curves.
In particular, since every complex K3 surface of Picard rank $\rho\geq19$
is rationally dominated by a Kummer surface,
these surfaces contain infinitely many rational curves.

In \cite{bt}, elliptic K3 surfaces $S\to\PP^1$ are studied. 
There, the authors define a {\em nt-multisection} to be a multisection 
$M$ of the fibration such that for a general point $b\in\PP^1$ there exist two points 
in the fiber
$p_b,p_{b'}\in S_b\cap M$ such that the divisor $p_b-p_{b'}$, considered as a point of
the Jacobian of $S_b$, is non-torsion.
Establishing infinitely many nt-multisections that are rational curves,
we find infinitely many rational curves on elliptic K3 surfaces of
Picard rank $\rho\leq19$, see \cite[Corollary 3.28]{bt}.
We note that K3 surfaces of Picard rank $\rho\geq5$ are automatically elliptic:
namely, in this case,
by the theory of integral quadratic forms, there exists an isotropic vector in $\Pic(S)$, 
which gives rise to an elliptic fibration.
Combining these results, we obtain the following

\begin{Theorem}[Bogomolov--Tschinkel]
  \label{bt elliptic thm}
  Let $S$ be a complex projective K3 surface that 
  \begin{enumerate}
    \item carries an elliptic fibration, or
    \item is a Kummer surface, or
    \item has Picard rank $\rho\geq 5$.
  \end{enumerate}
  Then, $S$ contains infinitely many rational curves.
\end{Theorem}

Moreover, in case the effective cone of a K3 surface
is not finitely generated, we find infinitely many rational curves
using Theorem \ref{one curve}.
Also, if the automorphism group is infinite, there are infinitely
many rational curves.
Combining these observations with the previous results, one can
show that there are infinitely many rational curves for K3 surfaces 
of Picard rank $\rho\geq4$, possibly with the exception of 
two Picard lattices of rank $4$.
We refer to \cite[Section 4]{bt} and \cite[Section 2]{bht} for details,
as well as to \cite[Example 4.8]{bt} for an example of a K3 surface
with $\rho=4$, where infinity of rational curves is currently still unknown.

On the other hand, a very general K3 surface in ${\cal M}_{2d}$ has Picard
rank $\rho=1$, does not carry an elliptic fibration, and has a finite 
automorphism group.
Thus, these are hard to come by, as they do not 
possess much geometric structure to work with.

\subsection{Reduction modulo $\mathbf{p}$}
In \cite{bht}, Bogomolov, Hassett and Tschinkel gave an approach to the case of
Picard rank $\rho=1$, which uses reduction modulo finite characteristic.
First, using degeneration techniques, they reduced to the number 
field case

\begin{Proposition}[Bogomolov--Hassett--Tschinkel]
  Bogomolov's conjecture \ref{infinity conjecture} holds for complex projective K3 
   surfaces if and only if it holds for K3 surfaces that are defined
   over number fields.
\end{Proposition}

Now, let $S$ be a K3 surface over some number field $K$.
Replacing $K$ by a finite extension, we may assume that all divisor classes
of $S_\CC$ are already defined over $K$.
Embedding $S$ into some projective space $\PP_K^N$, and taking the 
closure of its image
inside $\PP^N_{\OO_K}$, we get a model of $S$
over $\OO_K$.
After localizing at a finite set of places $P$ depending on $S$ and this
embedding, we obtain a smooth projective model 
${\cal S}\to\Spec\OO_{K,P}$, i.e., a smooth projective scheme over
$\OO_{K,P}$ with generic fiber ${\cal S}_K\iso S$.
In particular, for every prime ideal $\idealp$ of $\OO_{K,P}$, the reduction
${\cal S}_{\idealp}$ of $\cal S$ modulo $\idealp$ is a K3 surface over
the finite field $\OO_{K,P}/\idealp$.

The crucial observations and strategy of \cite{bht} are as follows:
let $(S,H)$ be a polarized K3 surface over $K$ with geometric Picard rank
 $\rho=1$, or, more generally,  $\rho$ odd.
If ${\rm char}(\OO_K/\idealp)\geq5$, then the Tate-conjecture holds for
${\cal S}_\idealp$ by Theorem \ref{tate conjecture}.
In particular, if we denote by ${\cal S}_{\overline{\idealp}}$ the base change of
${\cal S}_\idealp$ to the algebraic closure of ${\OO_K/\idealp}$,
then the Picard rank of ${\cal S}_{\overline{\idealp}}$ is even by
Theorem \ref{K3 even rank}.
On the other hand, the specialization map
$$
   \Pic(S)\,\iso\,\Pic({\cal S})\,\stackrel{{\rm sp}_\idealp}{\longrightarrow}\,
   \Pic({\cal S}_\idealp)
$$
is injective.
Since $\rho$ is odd, there exists for every prime $\idealp$
not lying over $2$ or $3$
an invertible
sheaf ${\cal L}_\idealp$ on ${\cal S}_{\overline{\idealp}}$ 
that does not lift to $S$.
We may assume ${\cal L}_\idealp$ to be effective, and then, by
Theorem \ref{one curve}, we find an effective divisor in
$|{\cal L}_\idealp|$ that is a sum of rational curves.
Since ${\cal L}_\idealp$ does not lift, there is at least one rational
curve $C_\idealp$ in this sum that does not lift to $S$ either. 
However, if $N_\idealp$ is a sufficiently large integer, then
$|N_\idealp H-C_\idealp|$ is effective, and by 
Theorem \ref{one curve}, there exist rational curves $R_{\idealp, i}$ on 
${\cal S}_{\overline{\idealp}}$ and positive
integers $n_i$ such that
\begin{equation}
  \label{stable curve}
  C_\idealp\,+\,\sum_i n_i R_{\idealp, i} \,\in\, | N_\idealp H |\,.
\end{equation}
This sum of rational curves can be represented by a stable map
of genus zero and so, defines a point of the moduli space
of stable maps
${\cal M}_0({\cal{S}}_{\overline{\idealp}}, N_\idealp H)$.
Next, we want this stable map to be {\em rigid}, i.e.,
the stable map allows at most infinitesimal deformations, i.e.,
the moduli space is zero-dimensional at this point.

The first problem is that rational curves can move on K3 surfaces in
positive characteristic (in which case we might not be able
to find a rigid representation).
But then, the K3 surface is uniruled, and in particular,
Artin-supersingular, see Section \ref{uni and super}.
By results of Bogomolov and Zarhin \cite{Bogomolov Zarhin}
(independently also obtained by Joshi and Rajan, but unpublished), 
we can always find infinitely many places $\idealp$ such that
${\cal S}_\idealp$ is not Artin-supersingular, which is sufficient for our
application.

Now, take of these infinitely many primes of non-supersingular reduction and
{\em suppose} (we comment on that below) that we can find a rigid
stable map representing (\ref{stable curve}).
We denote by $k$ the algebraic closure of the finite field
$\OO_K/\idealp$, let $W(k)$ be the Witt ring of $k$,
and base-change the family ${\cal S}\to \Spec\OO_{K,P}$ to
$W(k)$.
Then, dimension estimates of the relative formal moduli space 
${\cal M}_0({\cal{S}}, N_\idealp H) \to {\rm Spf}\: W(k)$
imply that our stable map to $S_{\overline{\idealp}}$ extends
to a stable map to the family ${\cal S}$ (here, rigidity 
is crucial).
Thus, the stable map lifts over a possibly ramified extension of
$W(k)$, and in particular, there exists a rational curve on $S_\CC$, 
whose reduction modulo $\idealp$ contains $C_\idealp$.
Thus, for infinitely many  $\idealp$ we get rational curves 
on $S_\CC$, and eventually obtain the following result
\cite{bht}:

\begin{Theorem}[Bogomolov--Hassett--Tschinkel]
  Let $S$ be a complex projective K3 surface with Picard group
  $\Pic(S)=\ZZ\cdot H$ such that $H^2=2$.
  Then, $S$ contains infinitely many rational curves.
\end{Theorem}

The main issue is the representation of (\ref{stable curve}) by
a {\em rigid} stable map, for otherwise it is not clear whether
one can lift this sum of rational curves to characteristic zero.

For degree $2$ and $\rho=1$, such a rigid representation exists
by exploiting the involution on K3 surfaces of degree $2$, see \cite{bht}.
In general, this difficulty was overcome in \cite{Li Liedtke} by introducing
{\em rigidifiers}: by definition, these are
ample and irreducible rational curves with at worst nodal 
singularities.
Then, every sum of rational curves can be represented by
a rigid stable map after adding sufficiently many rigidifiers to them.
Unfortunately, the surface $S_{\overline{\idealp}}$ may not contain rigidifiers.
However, surfaces containing rigidifiers are dense in the moduli space of
polarized K3 surfaces.
Using deformation techniques and rigidifiers, we obtained 
in \cite{Li Liedtke}

\begin{Theorem}[Li--Liedtke]
  Let $S$ be a complex projective K3 surface, whose Picard rank is odd.
  Then, $S$ contains infinitely many rational curves.
\end{Theorem}

More generally, the method of proof works whenever a K3 surface $S$ 
is defined over some field $K$, and we can find a DVR $R$ with quotient field 
$K$, as well as infinitely many primes $\idealp$ of $R$ 
such that the geometric Picard rank
of the reduction $S_\idealp$ is strictly larger than that of $S$.
For example, if $S$ is a complex projective K3 surface that cannot
be defined over a number field, then $S$ can be realized as generic
fiber of a non-isotrivial family ${\cal S}\to B$ over some positive dimensional 
base of characteristic zero.
Using results on the jumping of Picard ranks of K3 surfaces
in families from \cite{bkps} or \cite{oguiso},
we obtain

\begin{Theorem}
  Let $S$ be a complex projective K3 surface that cannot be defined
  over a number field. 
  Then, $S$ contains infinitely many rational curves.
\end{Theorem}

In view of these results and Theorem \ref{bt elliptic thm}, 
it remains to deal with K3 surfaces of Picard rank $\rho=2$ and
$\rho=4$ that are defined over number fields,
in order to establish Conjecture \ref{infinity conjecture} for all complex
projective K3 surfaces.
To apply the techniques of \cite{bht} and \cite{Li Liedtke}, we need jumping
of Picard ranks for infinitely places of non-supersingular reduction.
For example, such jumping results for certain classes of K3 surfaces with 
$\rho=2$ and $\rho=4$ over number fields were established in
\cite{Charles Picard}.

We end by giving a heuristic reason why we always expect 
to find infinitely many places with non-supersingular
reduction and jumping Picard rank (as in the case of odd rank), 
which would imply Conjecture \ref{infinity conjecture}.
However, in view of the results in \cite{maulik poonen}
and \cite[Theorem 1]{Charles Picard}, the situation
may be more subtle than expected.
In any case, here is our heuristic: 

The universal polarized K3 surface has Picard rank $\rho=1$.
All its (non-supersingular) specializations to surfaces over finite fields 
have a larger geometric Picard rank, and the extra invertible sheaves
extend (at least, formally) along divisors inside the moduli space. 
Also, these invertible sheaves must have unbounded intersection number 
with the polarization (otherwise some of them would lift to the universal K3
surface, which was excluded).
Thus, the moduli space of polarized K3 surfaces over the integers
is ``flooded''
by infinitely many divisors on which Picard ranks jump.
It is likely that given a K3 surface over a number field, infinitely many of its
non-supersingular reductions hit these divisors, establishing
the desired jumping behavior of Picard ranks.

\addtocontents{toc}{{\medskip \ }}

\end{document}